\documentclass[12pt]{article}
\pdfoutput=1

\newlength{\abstractwidth}
\usepackage{amsmath}
\usepackage{amsfonts}
\usepackage{amssymb}
\usepackage{latexsym}

\usepackage{epsf}
\usepackage{color}
\usepackage{graphicx}
\usepackage{tikz}
\usepackage{dsfont}

\usetikzlibrary{arrows,shapes,positioning}
\usetikzlibrary{decorations.markings}
\usepackage[rightcaption]{sidecap}
\tikzstyle arrowstyle=[scale=1]
\tikzstyle directed=[postaction={decorate,decoration={markings,
    mark=at position .65 with {\arrow[arrowstyle]{stealth}}}}]
\tikzstyle reverse directed=[postaction={decorate,decoration={markings,
    mark=at position .65 with {\arrowreversed[arrowstyle]{stealth};}}}]
\usetikzlibrary{positioning}

\renewcommand{\thefootnote}{\fnsymbol{footnote}}
\renewcommand{\thanks}[1]{\footnote{#1}}
\newcommand{\starttext}{
\setcounter{footnote}{0}
\renewcommand{\thefootnote}{\arabic{footnote}}}

\newcommand{\bea}{\begin{eqnarray}}
\newcommand{\eea}{\end{eqnarray}}
\newcommand{\be}{\begin{eqnarray}}
\newcommand{\ee}{\end{eqnarray}}
\newcommand{\bma}{\begin{matrix}}
\newcommand{\ema}{\end{matrix}}

\newtheorem{thm}{Theorem}[section]
\newtheorem{lem}[thm]{Lemma}
\newtheorem{prop}[thm]{Proposition}
\newtheorem{cor}[thm]{Corollary}
\newtheorem{conj}[thm]{Conjecture}

\def\cA{{\cal A}}
\def\cB{{\cal B}}
\def\cC{{\cal C}}

\def\cE{{\cal E}}

\def\cH{{\cal H}}

\def\cK{{\cal K}}
\def\cL{{\cal L}}
\def\cM{{\cal M}}

\def\cO{{\cal O}}

\def\cS{{\cal S}}
\def\cT{{\cal T}}

\def\mC{\mathfrak{C}}

\def\mI{\mathfrak{I}}

\def\mK{\mathfrak{K}}

\def\mN{\mathfrak{N}}

\def\mS{\mathfrak{S}}

\def\mm{\mathfrak{m}}

\def\mpp{\mathfrak{p}}

\def\ZZ{{\mathbb Z}}
\def\RR{{\mathbb R}}
\def\NN{{\mathbb N}}
\def\CC{{\mathbb C}}
\def\HH{{\mathbb H}}
\def\QQ{{\mathbb Q}}
\def\GG{{\mathbb G}}

\def\Re{{\rm Re \,}}
\def\Im{{\rm Im \,}}

\def\det{{\rm det \,}}

\def\half{{1\over 2}}
\def\p{\partial}

\def\G{\Gamma}

\def\f{\varphi}
\def\tet{\vartheta}
\def\ep{\varepsilon}

\def\pbz{\p _{\bar z}}

\def\th{{\rm th}}

\def\no{\nonumber}
\def\sm{\smallskip}

\begin{document}
\starttext
\setcounter{footnote}{0}

\begin{flushright}
2017 August 26 \\
2018 April 14 (revised)
\end{flushright}

\vskip 0.3in

\begin{center}

{\Large \bf Fourier series of modular graph functions }

\vskip 0.3in

{\large  Eric D'Hoker$^{(a)}$ and William Duke$^{(b)}$} 

\vskip 0.1in

 { \sl (a) Mani L. Bhaumik Institute for Theoretical Physics}\\
{\sl Department of Physics and Astronomy }

\vskip 0.1in

{\sl (b)  Department of Mathematics}\\
{\sl University of California, Los Angeles, CA 90095, USA} 

\vskip 0.05in

\vskip 0.07in

{\tt \small dhoker@physics.ucla.edu, wdduke@g.ucla.edu}

\end{center}

\begin{abstract}

Modular graph functions associate to a graph an $SL(2,\ZZ)$-invariant function on the upper half plane.  We obtain the Fourier series  of modular graph functions of arbitrary weight $w$ and two-loop order. The motivation for this work is to develop a deeper understanding of the origin of the algebraic identities between modular graph functions which have been discovered recently, and of the relation between the existence of these identities and the occurrence of  cusp forms.  We show that the constant Fourier mode, as a function of the modulus $\tau$,  consists of a Laurent polynomial  in $y = \pi \, \Im \tau$ of degree $(w,1-w)$, plus a contribution which decays exponentially as $y \to \infty$. The Laurent polynomial is a linear combination with rational coefficients of the top term $y^w$, and lower order terms $\zeta (2k+1) y^{w-2k-1}$ for $1\leq k \leq w-1$, as well as  terms $\zeta (2w-2\ell-3) \zeta (2\ell+1)y^{2-w}$ for $1 \leq \ell \leq w-3$.  The exponential contribution is a linear combination of exponentials of $y$ and incomplete $\Gamma$-functions whose coefficients are Laurent polynomials in $y$ with rational coefficients.  

\end{abstract}

\newpage


\newpage

\baselineskip=15pt
\setcounter{equation}{0}
\setcounter{footnote}{0}

\newpage

\section{Introduction and statement of the main result}
\setcounter{equation}{0}
\label{sec:1}

A modular graph function associates to a certain kind of  graph an $SL(2,\ZZ)$-invariant function of the upper half plane $\HH$. Modular graph functions naturally arise in the low energy expansion of closed string amplitudes and govern the contributions to this expansion at genus one \cite{Green:1999pv,Green:2008uj}. For one-loop graphs they are non-holomorphic Eisenstein series, while for two-loop graphs they were found to obey a system of differential equations \cite{D'Hoker:2015foa}. Connections between modular graph functions, multiple-zeta values and single-valued elliptic polylogarithms were put forward in \cite{DHoker:2015wxz}. Earlier relations between open and closed string amplitudes, multiple-zeta-functions, and polylogarithms were exhibited in \cite{Schlotterer:2012ny,Broedel:2013tta,Stieberger:2013wea}. The structure of the low energy expansion of genus-two string amplitudes, in as much as is known to date, may be found in \cite{DHoker:2014oxd} and references therein.  

\sm

A number of algebraic identities between modular graph functions with two or more loops were conjectured in~\cite{D'Hoker:2015foa} by matching their asymptotic expansions near the cusp $\tau \to i \infty$. The simplest of these identities were proven by direct summation of the Eisenstein series in \cite{Zagier:2014,D'Hoker:2015zfa}.  More complicated  identities and their various generalizations were obtained and proven by appealing to the rich system of differential equations they satisfy \cite{D'Hoker:2016a,DHoker:2016quv,Basu:2016xrt,Basu:2016kli}. In particular, it was shown in \cite{D'Hoker:2016a,DHoker:2016quv} that the existence of these identities may be traced back to the existence of identities between holomorphic modular forms.

\sm

Better understanding the origin of these algebraic identities is the main motivation for the present work, as will be explained further below.  Since the coefficients in the asymptotic expansion near the cusp of modular graph functions are known to  involve general multiple zeta-values \cite{DHoker:2015wxz,Zerbini}, the algebraic identities between them contain and generalize some of the well-known relations  between multiple zeta-values (see for example \cite{ZagierMZV, W, Zudilin, Brown:2013gia} and references therein) to the world of modular functions. A recent characterization of classes of modular graph functions, as single-valued projections of elliptic multiple zeta values introduced in \cite{enriquez}, may be found in  \cite{Brown:2017qwo,Brown2}.

\sm

A simple infinite family is given by modular graph functions depending  on $a_1,\dots a_\ell \in \NN$ with $\ell \geq 2$ and may be expressed for $\tau = \tau_1 + i \tau _2\in \HH$ with $\tau_1, \tau _2 \in \RR$ in the form,
\bea
\label{1a1}
C_{a_1,\cdots ,a_\ell} (\tau) = \sum _{{(m_r, n_r)\in \ZZ^2 \atop  r=1,2,\dots \ell}} ' \delta _{m,0} \, \delta _{n,0} \prod _{r=1}^\ell 
\left ( { \tau_2\over \pi |m_r  + n_r \tau |^2} \right )^{a_r}
\eea
where $m=m_1+m_2+\cdots m_\ell$, $n=n_1+n_2+\cdots n_\ell$ and the Kronecker $\delta$-symbols force $m,n$ to vanish. Throughout, a prime over the summation symbol indicates that the summation is restricted 
to disallow division by zero, in this case to $(m_r,n_r)\neq (0,0).$ The sum is absolutely convergent and $C_{a_1,\cdots ,a_\ell}  $ is smooth and modular, 
\bea
C_{a_1,\cdots ,a_\ell}(M \tau) =C_{a_1,\cdots ,a_\ell}(\tau) 
\eea
for each $M \in SL(2,\ZZ)$ acting on $\tau\in \HH$ by a linear fractional map.

\sm

 The associated graph $\G=\G_{a_1,\dots,a_\ell}$ is planar and may be realized by first taking a graph $\hat{\G}$ in the plane with 2 vertices connected to each other by $\ell$ edges, numbered $1,2, \dots , \ell$ and then adjoining $a_r-1$ distinct new bivalent vertices to the $r$th edge for each $r$. Thus $\hat{\G}$ has  $\ell$ edges, 2 vertices and $\ell$ faces, while $\G$ has  $w=a_1+\cdots + a_\ell$ edges, $w-1$ vertices and  $\ell$ faces, or equivalently $\ell-1$ loops.  We will review in the next section  how $C_{a_1,\cdots ,a_\ell} $, which is  said to be a  $(\ell-1)$-loop modular graph function of weight $w$, arises from $\G$. Its birth is as a Feynman graph in quantum field theory and string theory, and in that context the combination $p_r= m_r \tau + n_r$ represents the lattice momentum running through the edge $r$ taking values in the lattice $\Lambda= \ZZ + \ZZ \tau$. 

\sm

In the most simple  case of a one-loop graph we have $\ell=2$ and $C_{a_1,a_2}$ is given by a specialization of the classical  Kronecker--Eisenstein series,
since
\bea\label{eis1}
C_{a_1,a_2}(\tau) = E_w(\tau)=\sum _{(m,n) \in \ZZ^2} ' { \tau_2^w \over \pi^w |m+ n \tau|^{2w}}.
\eea
It is well known that $E_w(\tau)$ is an eigenfunction of the hyperbolic Laplacian,
\bea
\Delta E_w = w(w-1) E_w,
\eea
where  $\Delta$  is normalized by $\Delta = 4 \tau_2^2 \p_{\bar \tau} \p_\tau$. The effect of the Laplacian on $C_{a_1,\cdots ,a_\ell} (\tau)$ was obtained in  \cite{D'Hoker:2015foa} where it was also shown that, as a result, a number of identities between various kind of modular graph functions are forced to exist.
The simplest such identity was obtained in this way  in \cite{D'Hoker:2015foa} and states that,
\bea
\label{id1}
C_{1,1,1} (\tau)= E_3(\tau)+\zeta(3)
\eea
where $\zeta(s)$ is the Riemann zeta-function.  A  proof by direct summation of the Kronecker-Eisenstein series was given by  Zagier \cite{Zagier:2014}. In general, for $\ell =3 $ and for each odd value of $w$ there is precisely one linear combination of  those $C_{a_1, a_2, a_3}$ with weight $w$ that differs from  $E_w$ by a constant. For instance, for $w= 5, 7, 9$,
\bea
\label{id2}
30C_{2,2,1} (\tau) & = & 12 E_5 (\tau)  + \zeta (5) 
\no \\
252C_{3,3,1} (\tau)+252C_{3,2,2}(\tau) & = & 108E_7(\tau)+ \zeta (7) 
\no \\
2160 C_{4,4,1}(\tau)  + 4320 C_{4,3,2}(\tau) + 960 C_{3,3,3}(\tau) & = & 960 E_9(\tau) +\zeta (9) 
\eea
Note that, if we assign ``weight $s$" to $\zeta (s)$, then the above identities are all homogeneous of their respective weights. The case of even weight  is more difficult. Here the identities involve the functions $C$ with $\ell >3.$ For instance we have, 
\bea
\label{id3}
C_{1,1,1,1}(\tau)=24C_{2,1,1}(\tau)-18E_{4}(\tau)+3E_2^2(\tau)
\eea
This identity was conjectured in  \cite{D'Hoker:2015foa} and proven in \cite{D'Hoker:2015zfa}.

\sm

Identities such as (\ref{id1}), (\ref{id2}) and (\ref{id3})  are reminiscent of various identities from the classical theory of modular forms. Consider the classical holomorphic Eisenstein series of even modular weight\footnote{The {\sl modular weight} is in general distinct from the weight $w$ of the modular functions $C$ defined earlier.}  $w\geq 4$,
\bea
\GG_{w}(\tau)= \frac{(w-1)!}{2 (2\pi i)^w} \sum'_{(m,n)\in \ZZ^2}{ 1 \over (m+n\tau)^w} 
\eea
We have  the Fourier expansion,
\bea
\GG_{w}(\tau)=-\frac{B_{w}}{2w}+ \sum_{k \geq 1} \sigma_{w-1}(k) q^k,
\eea
where $B_{w}$ is the Bernoulli number, $q=e^{2\pi i \tau}$ and $\sigma_{s}(n)= \sum_{d|n} d^{s}$ is the divisor sum. The following identity is forced  by the fact that $\GG_8$ and $\GG_4^2$ are in the same one-dimensional space,
 \bea
120\,\GG_4^2 =\GG_8
\eea
and implies the following arithmetic identity,
\bea
\sigma_7(n)= \sigma_3(n)+120 \sum_{m=1}^{n-1} \sigma_3(n-m)\sigma_3(m).
\eea
At higher weight such identities will involve holomorphic cusp forms, a phenomenon which arises first at weight 12 where we have, for example,
\bea
65\,\GG_{12}(\tau)-174132\,\GG_6^2(\tau) = 756\, \Delta(\tau),
\eea
where $\Delta(\tau)=q\prod_{m}(1-q^m)^{24}$ is the famous cusp form of modular weight 12. Knowledge of the dimensions of holomorphic modular forms for given modular weight, together with their Fourier series expansion, provides all the information needed to establish all such identities in the holomorphic case. 

\sm

Returning to non-holomorphic modular graph functions, it is one of the ultimate goals of this and subsequent  work to understand the relation between identities amongst modular graphs functions, such as (\ref{id1}), (\ref{id2}) and (\ref{id3}), their Fourier series expansion, and the existence of non-holomorphic cusp forms. The Fourier expansion  of the non-holomorphic Eisenstein series $E_w$ for integer $w \geq 2$ is well-known and is given by,
\bea
\label{2d5}
E_w(\tau) & = & 
- { B_{2w} \over (2w)!} (-4y)^w + { 4 (2w-3)! \over (w-2)! \, (w-1)!} \, \zeta (2w-1) \, (4y)^{1-w} 
\no \\ &&
+ { 2 \over (w-1)!} \sum _{k=1} ^\infty  k^{w-1} \sigma _{1-2w}(k) \, \Big ( q^k + \bar q^k \Big )
P_w(4 k y).
\eea
where  we set $y = \pi \tau_2$, and $P_w(x) $ is the polynomial in $1/x$ defined by,
\bea
\label{Pw}
P_w (x) = \sum _{m=0}^{w-1} { (w+m-1)! \over m! \, (w-m-1)! ~ x^m}
\eea 
For example,  we have, 
\bea
\label{FE1}
E_3(\tau)  = 
 \tfrac{2}{945} y^3 + \tfrac{3}{4} \zeta (5) \, y^{-2} 
+  \sum _{k=1} ^\infty \sigma _{3}(k)\Big(\tfrac{1}{k}+\tfrac{3}{2k^2y}+\tfrac{3}{4k^3y^2}\Big) \, \Big ( q^k + \bar q^k \Big ).
\eea
We see that the constant Fourier mode of $E_w$ is a Laurent polynomial with only two terms. 

\sm

In the present paper, we shall study the Fourier series expansion of two-loop modular graphs functions $C_{a_1, a_2, a_3}$. By combining (\ref{id1}) and (\ref{id2}) with (\ref{FE1}) we readily get the Fourier expansions of $C_{1,1,1}$ and $C_{2,2,1}$. For functions of higher weight, the Laurent series have been evaluated only in special cases, tabulated below, 
\bea
\label{spec}
C_{2,1,1} (\tau)
& = &   \frac{2   y^4 }{14175}  +  \frac{ \zeta(3) y}{45}  + \frac{5  \zeta(5)}{12 y}  -\frac{ \zeta(3)^2}{4 y^2} +  \frac{9 \zeta(7)}{16  y^3} + \cO(e^{-2 \pi \tau_2}) 
\no \\
C_{3,1,1} (\tau) 
& = &
\frac{2  y^5}{155925}  + \frac{2   \zeta (3) y^2 }{945}  -\frac{\zeta (5)}{180}
+\frac{7 \zeta (7)}{16 y^2} -\frac{\zeta (3)  \zeta (5)}{2 y^3} +\frac{43  \zeta (9)}{64 y^4} + \cO(e^{-2 \pi \tau_2}) 
\no \\
C_{4,1,1} (\tau) &=&\frac{808 y^6}{638512875}+\frac{ \zeta (3)  y^3}{4725}-\frac{ \zeta (5)  y}{1890} +\frac{\zeta (7)}{720 y} +\frac{23 \zeta (9)}{64y^3} 
  \no\\&&
-\frac{\zeta (5)^2+30 \zeta(3) \zeta (7)}{64 y^4}+\frac{167 \zeta (11)}{256 y^5} + \cO(e^{-2 \pi \tau_2}) 
  \no \\
C_{3,2,1}(\tau) &=&\frac{43 y^6}{58046625}+\frac{y \zeta (5)}{630}+\frac{\zeta (7)}{144 y}+\frac{7 \zeta (9)}{64 y^3}-\frac{17 \zeta (5)^2}{64y^4}+\frac{99 \zeta (11)}{256 y^5} + \cO(e^{-2 \pi \tau_2})
  \no \\
C_{2,2,2} (\tau) &=&\frac{38 y^6}{91216125}+\frac{\zeta (7)}{24 y}-\frac{7 \zeta (9)}{16 y^3}+\frac{15 \zeta (5)^2}{16y^4}-\frac{81 \zeta (11)}{128 y^5} + \cO(e^{-2 \pi \tau_2}) 
 \eea
Along different lines, a systematic algorithm was developed to evaluate modular graph functions which have at most four vertices (including bivalent vertices) in \cite{Zerbini}. 

\sm

Except in such cases where we may apply a known identity, it is an open problem to determine the Fourier coefficients of $C_{a_1,\cdots ,a_\ell} (\tau)$ for a fixed $a_1,\cdots ,a_\ell$ when $\ell >2$. Its expansion has the form,
\bea
\label{CF}
C_{a_1,\cdots ,a_\ell} (\tau) = \sum _{k=- \infty } ^\infty C_{a_1, \cdots, a_\ell}^{(k)}  ( \tau_2) \, e^{2 \pi i k \tau_1},
\eea
If we can determine $C_{a_1, \cdots, a_\ell}^{(k)}  ( \tau_2)$, independently of any knowledge of special identities between modular graph forms, then  identities such as (\ref{id1}), (\ref{id2}) and (\ref{id3}) and their generalizations should emerge from identities between the Fourier coefficients of these modular graph functions.
Understanding the structure of the Fourier series thus appears key to understanding the structural mechanism behind the existence of the identities between modular graph functions, and their relation with the existence of 
cusp forms. Clearly, one obstruction to the existence of identities is the presence of cusp forms. Therefore, an urgent question is whether we can find combinations of the $C$-functions of a fixed weight having zero constant mode in their  Fourier expansions. Finally, it would be interesting to find out whether such cusp forms play a natural role  in string theory.

\sm

As a  step toward answering  such questions,  in this paper we will compute rather explicitly the constant mode $C_{a_1, a_2, a_3}^{(0)}  ( \tau_2)$ in the  Fourier expansion of $C_{a_1,a_2,a_3}$.   The main results of this paper may be summarized by three Theorems, and one conjectured decomposition formula.

\medskip
\noindent
\begin{thm}\label{t1}
 The constant Fourier mode for $C_{a_1,a_2,a_3}$ with fixed $a_1,a_2,a_3$ is given by,
\bea
\label{CE}
C_{a_1, a_2, a_3}^{(0)}  ( \tau_2) = \cL(\tau_2) + \mathcal{E}(\tau_2)
\eea
where $\cL(\tau_2) $ is a Laurent polynomial in $\tau_2$ of degree $(w,1-w)$ and $ \mathcal{E}(\tau_2) $ is exponentially decaying as $\tau \to i \infty$. The Laurent polynomial is given by, 
\bea
\label{CL}
\cL(\tau_2) = c_w (-4 \pi \tau_2)^w + \sum _{k=1}^{w-1} c_{w-2k-1} {\zeta (2k+1) \over ( 4 \pi \tau_2)^{2k+1-w}}  
+{c_{2-w}  \over (4 \pi \tau_2)^{w-2}}
\eea 
\begin{description}
\item (a) The coefficient $c_w$ is a rational number given by, 
\bea
c_w = \sum _{k=0} ^{ a_2}   {   B_{2k} B_{2w-2k} \over (2k)! \, (2w-2k)! } \, { \Gamma (2a_2+2a_3-2k) \over \Gamma (2a_3) \Gamma (2a_2-2k+1)}
+ (a_2 \leftrightarrow a_3)
\eea
\item (b) The coefficients $c_{w-2k-1}$ for $1 \leq k \leq w-1$ are rational numbers given by, 
\bea
c_{w-2k-1}   & = &  
 { 2 B_{2w-2k-2}    \over  (2w-2k-2)! } \sum _{\alpha =0}^{a_1-1} \sum _{\beta =0} ^{a_1-1-\alpha}  
 (-)^{a_1+a_3+\beta+1} \, \theta \left ( a_3+\left [{a_2+\beta \over 2}  \right ] - w+k+1 \right ) 
 \no \\ && 
\times g_{a_1,a_2} (\alpha, \beta) 
\binom{2k-2a_1+\alpha + \beta +1}{a_2+\alpha -1}
+ \hbox{5 permutations of } a_1, a_2, a_3
\no
\eea 
where $\theta(x)$ is the step function defined to equal 1 when $x\geq 0$ and to vanish otherwise. The function $g_{a_1, a_2}(\alpha, \beta)$ is integer-valued and given by,
\bea
\label{G3}
g_{a_1,a_2} (\alpha, \beta) = 
(-)^{a_1} 
\left ( \bma 2a_1-2-\alpha -\beta \cr a_1-1 \cr \ema \right )
\left ( \bma a_2+ \alpha -1 \cr a_2-1 \cr \ema \right )
\left ( \bma a_2+ \beta -1 \cr a_2-1 \cr \ema \right )
\eea
\item (c) The coefficient $c_{2-w}$ is given by the following expression, 
\bea
\label{C2W}
c_{2-w} = c^0_{2-w} (-)^w \zeta (2w-2)
 + 2 \sum_{\sigma \in \mS_3} Z(a_{\sigma(1)}, a_{\sigma(2)}, a_{\sigma(3)} )
\eea
Here, $c^0_{2-w}$ is an integer given by, 
\bea
\label{C2W0}
c^0_{2-w} = 
 \sum _{\alpha =0} ^{a_1-1} \sum _{\beta =0}^{a_1-1-\alpha} (-)^{a_2+\beta} g_{a_1,a_2}(\alpha, \beta) 
 \binom{2a_2+2a_3+\alpha + \beta}{a_2+\alpha -1}
\eea
while $Z(a_1,a_2,a_3)$ is a linear combination with integer coefficients of depth-two multiple zeta-functions of total weight $2w-2$ given by,
\bea
\label{ZZ}
Z(a_1,a_2,a_3)
& = & 
\sum_{k=1}^{a_1} \sum _{\ell=1}^{a_1} 
\binom{a_1+a_2-k-1}{a_2-1} \binom{a_1+a_2-\ell-1}{a_2-1} \binom{k+\ell-2}{k-1}
\no \\ && \hskip 0.4in \times 
 \binom{2w-k-\ell-2}{w-k-1} \, \zeta (2w-k-\ell-1, k+\ell-1)
\eea
\end{description}
\end{thm}
The normalization of the double $\zeta$-functions is as follows,
\bea
\label{zeta2}
\zeta (a,b)= \sum _{m,n=1} ^\infty \frac{1}{(m+n)^a n^b}
\eea

\medskip
\noindent
\begin{thm}\label{t2}
 The coefficient $c_{2-w}$ is a linear combination, with integer coefficients, of products of two odd  zeta-values whose weights add up to $2w-2$, 
 \bea
 \label{thm2}
c_{2-w} = \sum _{k=1}^{ w-2}  \half \gamma _k \, \zeta (2k+1) \zeta (2w-2k-3) 
\eea
with rational coefficients $\gamma _k \in \QQ$.
\end{thm}

\medskip
\noindent
\begin{thm}\label{t3}
The general structure of the exponential part (\ref{CE}) is given as follows, 
\bea
\mathcal{E}(\tau_2)= \sum _{s=1}^{w-2} \sum _{n=1} ^\infty \left ( \left (f_+(s;n) (4\pi \tau_2)^s + { f_-(s;n) \over (4 \pi \tau_2)^{s-1}} \right ) \, {\rm Ei }(4 \pi n \tau_2) + e^{-4 \pi n\tau_2} \sum _{m=2-w} ^{w-3} { f(s;n,m) \over (4\pi \tau_2)^m} \right ) \quad
\eea
where ${\rm Ei}$ is the incomplete $\Gamma$-function and  the coefficients $f_\pm(s,n)$,  and $f(s,n,m)$ are rational numbers.
\end{thm}

\medskip
\noindent
\begin{conj}
\label{conj1} (Decomposition Formula)
The coefficients $\gamma_k$ entering the decomposition of $c_{2-w}$ in formula (\ref{thm2}) of Theorem \ref{t2} are given by the following expression, 
\bea
\label{sumL3}
\gamma _k & = &  2 Z_k (a_1, a_2, a_3) \theta (a_1-1-k) - Z_0(a_1, a_2, a_3)  
\no \\ && 
+  \sum _{\alpha =1}^{a_1-1} Z_\alpha (a_1, a_2, a_3) \sum _{n=0}^{2\alpha -1} E_n(0) \binom{2k}{2\alpha -n} \binom{2w-2\alpha+n-4}{n}
\no \\ &&
+ \hbox{ 5 permutations of } a_1, a_2, a_3
\eea
and are integers.  Here, $E_n(x)$ are the Euler polynomials and the integer-valued function $Z_\alpha(a_1, a_2, a_3)$ is given by the following sum,
{\small \bea
\label{sumL4}
Z_\alpha  (a_1, a_2, a_3) = 
\sum _{k=k_-}^{k_+} \binom{a_1+a_3-k-1}{a_3-1} \binom{a_1+a_3-2\alpha +k -3}{a_3-1}
\binom{2\alpha}{k-1} \binom{2w-2\alpha-4}{w-k-1}
\eea}
with $k_+ = \min(a_1, 2\alpha+1)$ and $k_- = \max(1, 2 \alpha +2-a_1)$.
\end{conj}

\sm

 To obtain the decomposition formula, we make use of a conjectured relation (conjecture \ref{conj2}), which we have verified extensively using Maple calculations, but for which we have no analytical proof. Therefore, a full proof of the decomposition formula remains outstanding.

\sm

The explicit formulas for the Laurent polynomial part $\cL$ of $C_{a_1, a_2, a_3}(\tau)$, obtained in Theorems \ref{t1} and \ref{t2}, and the Decomposition Formula of \ref{conj1}, completely reproduce  the Laurent polynomials of (\ref{spec})  which have been evaluated earlier in the literature.

\section{Modular graph functions}
\setcounter{equation}{0}
\label{sec:2}

Before turning to the proof of the Theorems, and the derivation of the Decomposition Formula, we will briefly review the general definition of a modular graph function as it comes from string theory and show that the function $C_{a_1,\cdots , a_\ell}(\tau)$, which was defined in (\ref{1a1}),  is one such modular graph function.

\sm

Modular graph functions arise  as follows. The torus  $\Sigma$ with modulus $\tau $ may be represented in the complex plane by the quotient $\Sigma = \CC / \Lambda$ for the lattice $\Lambda = \ZZ + \ZZ \tau$. We choose local complex coordinates $(z, \bar z) $ on $\Sigma$ in which the metric  is given by $|dz|^2/\tau_2 $.  The volume form  of this metric $d\mu (z) =  i d z \wedge d \bar z/ (2 \tau_2)$ has  unit area and the Dirac $\delta$-function $\delta (z-w)$ is normalized by $\int _\Sigma d\mu(z) \delta (z-w) = 1$. The scalar  Green function $G(z-w|\tau)$  is defined by,
\bea
\label{2b1}
\tau_2 \pbz \p_z \, G(z-w|\tau) = -\pi \delta (z-w) + \pi
\eea
along with the normalization condition,
\bea
\label{2b2}
\int _\Sigma d\mu (z) \, G(z-w|\tau) = 0
\eea

To a graph $\Gamma$ with $v\geq 2$ vertices and  $w$ edges we associate $v$ points $z_i$ on the torus $ \Sigma$ labelled by the index $i=1,\cdots, v$. We denote by $\nu_{ij}$ the number of edges connecting the pair of  vertices $i,j =1, \cdots, v$. The number $\nu_{ij}$ is allowed to be a positive or zero integer for any pair of distinct vertices $i,j$. We set $\nu_{ii}=0$ for all vertices $i=1,\cdots , v$, and thus restrict the type of graphs on which we can define modular graph functions.\footnote{An immediate justification for this restriction on the graphs is that without it divergent contributions  involving the Green function at coincident points $G(0|\tau)$ would arise. In quantum field theory such graphs do arise in un-renormalized correlation functions, but are eliminated by the process of renormalization.}  The total number of edges is $w = \sum _{1 \leq i < j \leq v} \nu_{ij}.$  The modular graph function $\cC_\Gamma (\tau) $ is defined in terms of  absolutely convergent integrals over the torus by
\bea
\label{2b4}
\cC_\Gamma (\tau) = \left ( \prod _{k=1}^v \int _\Sigma d \mu (z_k) \right ) \prod _{ 1 \leq i < j \leq v} G(z_i-z_j|\tau)^{\nu_{ij}}
\eea
$\cC_\Gamma (\tau)$ is  clearly modular.

\sm

We may assume that $\G$ is connected and remains connected after the removal of a vertex (and its adjoining edges), since otherwise $C_\G=C_{\G_1}C_{\G_2}$ for subgraphs $\G_1$ and $\G_2.$ We may also assume that $\G$  remains connected when any single edge is omitted, for otherwise $C_\G=0$.
In particular we may assume that $\G$ contains no vertices with valence 1.  Suppose that $\G$ has $\hat{v}$ vertices of valence at least 3. As in the case of the graphs $\G_{a_1,\cdots ,a_\ell}$, which were introduced in the third paragraph of section \ref{sec:1},  it is convenient to build up $\G$ using an auxiliary graph $\hat{\G}$ having $\hat{v}$ 
vertices with the same valences $\geq 3$, no bivalent vertices and $\ell$ edges.  Note that, unlike $\G$,  $\hat{\G}$ can be a single edge with no vertices. This happens when $\hat{v}=0.$
To recover $\G$  we adjoin $a_r-1$ bivalent vertices to the $r$th edge of $\hat{\G}$   for $r=1,\dots \ell.$
 If $\hat{v}>0$ consider the $\hat{v}\times \ell$ incidence matrix of the graph $\hat{\G}$, when it is given some orientation.
The entry  $a_{i r}$ of the incidence matrix equals $\pm 1$ if edge $r$ starts or ends on vertex~$i$ (the sign is determined by the choice of orientation  through the graph), and equals 0 otherwise. 

\begin{prop}\label{p1}
Under the assumptions and notation introduced above, we have,
\bea
\label{2c6}
\cC_\Gamma (\tau) = \sum _{{ (m_r,n_r) \in \ZZ^2 \atop r=1,\cdots,\ell}} '  \, \prod _{r =1} ^{\ell} \left ( { \tau_2 \over \pi  |m_r + n_r \tau|^{2} }\right )^{a_r} \, 
\prod _{i=1}^{\hat v} \delta \left ( \sum _{r =1}^{\ell} a_{i r} m_r \right )
\delta \left ( \sum _{r =1}^{\ell} a _{i r} n_r \right )
\eea
When $\hat{v}=0$ we must evaluate the second product to 1.
\end{prop}

{\it Proof:}
The Green function $G (z|\tau)$, defined in (\ref{2b1}) and (\ref{2b2}), is given by a Fourier sum on the torus,
parametrized in terms of real coordinates $x,y$ by $z = x+ y \tau$ and $x,y \in \RR/\ZZ$, 
\bea
G(z|\tau) = \sum _{(m,n) \in \ZZ^2} ' { \tau_2  \over \pi |m + n \tau|^2} \, e^{2 \pi i (my-nx) }
\eea
Bivalent vertices play a special role, as they produce a convolution of concatenated Green functions. We parametrize their effect by introducing the functions $G_a (z|\tau)$, defined recursively in the index $a$ by setting $G_1 (z|\tau) = G(z|\tau)$ for $a=1$ and,
\bea
G_a (z|\tau) = \int _\Sigma d \mu(w) \, G(z-w|\tau) \, G_{a-1}(w |\tau)
\eea
for $a \geq 2$.  The Fourier series for $G_a$ on the torus  is readily obtained using  $d\mu(z) = dx \wedge dy$,
 \bea
 \label{2c3}
G_a(z|\tau) = \sum _{(m,n) \in \ZZ^2} ' { \tau_2^a  \over \pi^a |m + n \tau|^{2a}} \, e^{2 \pi i (my-nx) }.
\eea
Applying this to (\ref{2b4})   and carrying out the integrals over the $\hat{v}$ vertex positions $z_i$ corresponding to vertices of valence $\geq 3$, we obtain (\ref{2c6}).

\begin{cor} For $C_{a_1,\cdots ,a_\ell}$ and $\G_{a_1,\cdots ,a_\ell}$ defined in and below (\ref{1a1}) we have,
\bea
C_{\G_{a_1,\cdots ,a_\ell}}=C_{a_1,\cdots ,a_\ell}
\eea
\end{cor}

\section{Fourier series of two-loop modular graph functions}
\setcounter{equation}{0}
\label{sec:3}

Now we turn to the proof of the Theorems. In this section we shall introduce a Mellin-transform formulation of two-loop modular graph functions $C_{a_1, a_2, a_3}$  and a partial Poisson resummation to obtain the Fourier series expansion of (\ref{CF}). The method naturally generalizes to the case of higher modular graph functions, but we shall treat here only the case of two-loop modular graph functions of arbitrary weight $w=a_1+a_2+a_3$.

\subsection{Mellin-transform representation}

We begin with the elementary integral representation, 
\bea
{ \tau_2^{a_r}  \over \pi^{a_r}  |m_r + n_r \tau|^{2a_r}}
= 
 \int _0 ^\infty dt_r \, \frac{t_r ^{a_r-1}}{\Gamma (a_r)}   \exp \left \{ - { \pi \over \tau_2} t_r |m_r + n_r \tau|^2 \right \} 
 \eea
Collecting the sum over the product of three such factors, and taking care of omitting the zero mode from the summation over each edge, we find, 
\bea
\label{CS}
C_{a_1, a_2, a_3} (\tau) & = & 
\left ( \prod _{r=1}^3 \int _0^\infty dt_r \, \frac{t_r ^{a_r-1}}{\Gamma (a_r)}  \right ) \, S(t_1, t_2, t_3 |\tau)
\no \\
S(t_1,t_2,t_3|\tau ) & = & \sum _{{(m_r, n_r) \in \ZZ^2 \atop r=1,2,3}}  \delta _{m,0} \, \delta _{n,0}
\prod_{r=1}^3 \left (  \exp \left \{ - { \pi \over \tau_2} t_r |m_r + n_r \tau|^2 \right \} - \delta _{m_r,0} \delta _{n_r,0} \right ) \quad
\eea
The function $S(t_1, t_2, t_3|\tau)$ is  invariant under $SL(2,\ZZ)$ acting on $\tau$, as well as under permutations of the $t_r$. Next, we decompose $S$ by expanding the triple product into a sum of eight terms. The three contributions for which two pairs $(m_r, n_r)$ are set to zero must also have the third pair equal to zero in view of overall momentum conservation, and therefore combine with the terms in which all three pairs are zero. The result is as follows,
\bea
\label{SAB}
S(t_1, t_2, t_3|\tau ) = A(t_1, t_2, t_3|\tau ) - B(t_1+t_2|\tau ) - B(t_2+t_3|\tau ) -B(t_3+t_1|\tau ) +2 
\eea
The last term arises from the contribution with all pairs $(m_r,n_r)$ equal to $(0,0)$, and the functions $A$ and $B$ are given by,
\bea
A(t_1,t_2,t_3|\tau ) & = & \sum _{{(m_r, n_r) \in \ZZ^2 \atop r=1,2,3}}  \delta _{m,0} \, \delta _{n,0} \, 
\exp \left \{ - { \pi \over \tau_2} \sum _{r=1}^3 t_r |m_r + n_r \tau|^2 \right \} 
\no \\
B(t|\tau ) & = & \sum _{(m_1,n_1) \in \ZZ^2} \exp \left \{ - { \pi \over \tau_2}  t |m_1 + n_1 \tau|^2 \right \}
\eea
Note that the summation over the three pairs $(m_r,n_r)$ in $A$ includes the contribution from all the zero pairs, and is constrained only by the requirement that their sum $(m,n)$ vanishes. The summation over pairs $(m_1,n_1)$ in $B$ is unconstrained. 

\sm

It will be convenient to solve the constraint $m=n=0$  in the summation which defines the function $A$ by setting $m_3=-m_1-m_2$ and $n_3=-n_1-n_2$, with $m_1,m_2, n_1, n_2$ taking values in $ \ZZ$ unconstrained. Furthermore, we introduce the matrix notation, 
\bea
M = \left ( \bma m_1 \cr m_2 \cr \ema \right )
\hskip 0.6in
N = \left ( \bma n_1 \cr n_2 \cr \ema \right )
\hskip 0.6in
T = \left ( \bma t_1+t_3 & t_3 \cr t_3 & t_2 + t_3 \cr \ema \right )
\eea
The function $A$ then takes the form, 
\bea
A(t_1,t_2,t_3|\tau ) = \sum _{M,N \in \ZZ^2}  
\exp \left \{ - { \pi \over \tau_2} (M+ \tau N) ^\dagger T (M+ \tau N) \right \} 
\eea
One may think of this expression as defining a $\tet$-function.

\subsection{Partial Poisson resummation}

To compute the Fourier series of $C_{a_1, a_2, a_3}(\tau)$ as a function of $\tau_1$, we  perform a Poisson resummation on the sum in the expression for the function $A$ on the matrix $M$, but {\sl not} on~$N$. To do so, we  evaluate the Fourier transform of the $M$-dependent part as follows,
\bea
\int _{\RR^2} d^2M  \, e^{- 2 \pi i M^t X} e^{ - \pi (M+ \tau_1 N)^t T (M+ \tau_1 N) /\tau_2} 
= { \tau_2 \over (\det T)^\half } \, e^{ 2 \pi i \tau_1 N^t X - \pi \tau_2 X^t T^{-1} X }
\eea
It will be convenient to express the inverse of $T$ as follows, 
\bea
T^{-1} =  \frac{1}{\det T} \,  \ep^t \, T \, \ep 
\hskip 1in 
\ep =  \left ( \bma 0 & 1 \cr -1 & 0 \cr \ema \right )
\eea
Upon the change of summation variables $M \to - \ep M$,  the Fourier series takes the form, 
\bea
\label{At}
A(t_1, t_2, t_3 |\tau )  =   { \tau_2 \over (\det T)^\half } \sum _{M,N \in \ZZ^2} \,  e^{  2 \pi i M^t \ep N \tau_1 }  
 \exp \left \{    - \pi  \tau_2   {M^t T M \over \det T}   - \pi \tau_2 N^t  T  N   \right \}
\eea
To obtain the Fourier series of the function $C_{a_1, a_2, a_3}(\tau)$ we shall need to integrate $S$ over $t_1, t_2, t_3$, which requires combining the contributions of $A$ to the integral with those from $B$. 
To simplify this recombination, we perform a Poisson resummation in $m_1$ of $B(t)$,
\bea
\label{Bt}
B(t |\tau ) = \sqrt{{ \tau_2 \over t}} \sum _{m_1,n_1 \in \ZZ} e^{ 2 \pi i m_1 n_1 \tau_1} \, e^{- \pi \tau_2 m_1^2/t  - \pi \tau_2  n_1^2 t}
\eea
which exhibits the Fourier series  in $\tau_1$ of $B$.

\subsection{Fourier series expansions of $A$ and $B$}

The Fourier modes $S_k (t_1, t_2, t_3 |\tau_2)$ of  $S(t_1, t_2, t_3 |\tau)$ as a function of $\tau_1$ are given   by, 
\bea
S(t_1, t_2, t_3 |\tau)  =  \sum _{ k \in \ZZ} e^{2 \pi i k \tau_1} S_k (t_1, t_2, t_3|\tau_2)
\eea
The Fourier modes $A_k (t_1, t_2, t_3 |\tau_2)$ of $A(t_1, t_2, t_3 |\tau)$, and the Fourier modes $B_k (t |\tau_2)$ of  $B(t|\tau)$ as functions of $\tau_1$ are defined analogously. They are related to one another by,
\bea
\label{SABk}
S_k(t_1, t_2, t_3|\tau_2) & = & A_k(t_1, t_2, t_3|\tau_2) - B_k(t_1+t_2|\tau_2) 
\no \\ &&
- B_k(t_2+t_3|\tau_2) -B_k(t_3+t_1|\tau_2) + 2 \delta _{k,0}
\eea
The expressions for the Fourier modes are obtained from (\ref{At}) and (\ref{Bt}) and are given by,
\bea
\label{AkBk}
A_k (t_1, t_2, t_3 |\tau_2 ) & = &  { \tau_2 \over (\det T)^\half } \sum _{M,N \in \ZZ^2 }  \delta _{ M^t \ep N, k}    \,
 \exp \left \{    - \pi  \tau_2  {M^t T M \over \det T}   - \pi \tau_2  N^t  T  N   \right \}
 \no \\
B_k(t |\tau_2 ) & = & \sqrt{{ \tau_2 \over t}} \sum _{m,n \in \ZZ}  \delta _{mn,k}  \, e^{- \pi \tau_2 m^2/t  - \pi \tau_2  n^2 t} 
\eea
The Fourier modes $B_k$ for $k \not= 0$ are exponentially decaying as $t \to \infty$. Thus, the term $B_k(t_1+t_2|\tau)$ decays exponentially as $t_1 \to \infty$ or $t_2 \to \infty$ or both, but not when $t _3 \to \infty$. We will show in the subsequent subsection that uniform exponential decay is recovered  upon combining the contributions of the Fourier modes $A_k$ and $B_k$ into the modes $S_k$ given by (\ref{SAB}). With exponential decay secured, the Fourier modes of $C_{a_1, a_2, a_3}(\tau)$, expressed  in the notation of (\ref{CF}) with the help of (\ref{CS}), are then obtained by the following integrals,  
\bea
C_{a_1, a_2, a_3} ^{(k)} (\tau_2)=
\left ( \prod _{r=1}^3 \int _0^\infty dt_r \, \frac{t_r ^{a_r-1}}{\Gamma (a_r)}  \right ) \, S_k(t_1, t_2, t_3 |\tau)
\eea
which are absolutely convergent for large $t_r$, and may be analytically continued in $a_r$ for small $t_r$ if necessary.

\subsection{Partitioning the sum over $N$}

To expose uniform exponential decay in $t_1, t_2, t_3$, we partition the summation over $N \in \ZZ^2$, 
\bea
\label{part}
\ZZ^2 & = & \mN^{(0)} \cup \mN^{(1)} \cup \mN^{(2)} \cup \mN^{(3)}  \cup \mN^{(4)}
\eea
into the following (disjoint) parts,
\bea
\mN^{(0)} & = & \{ (0,0)\} 
\no \\
\mN^{(i)} & = & \{ (n_1, n_2) \hbox{ such that } n_1, n_2 \in \ZZ, n_i=0, n_j\not=0 \hbox{ for } j \not= i \}
\hskip 0.4in i,j=1,2,3
\no \\
\mN^{(4)}  & = & \{ (n_1, n_2) \hbox{ such that }  n_1, n_2 \in \ZZ, n_1, n_2, n_3 \not=0 \}
\eea
where we enforce the constraint $n=n_1+n_2+n_3=0$ throughout.  The quadratic form $N^t T N$ vanishes for $N\in \mN^{(0)}$; is uniformly non-degenerate in $t_1, t_2, t_3$ for  $N \in \mN^{(4)}$; while for $i=1,2,3$ and $N \in \mN^{(i)}$ the quadratic form decays exponentially, but non-uniformly,  in all directions except $t_i \to \infty$ where it remains bounded. Therefore, the above partition accurately governs the asymptotic  behavior as $t_1, t_2, t_3 \to \infty$.  We arrange  the contributions to $S_k$ arising from the  partitions $\mN^{(i)}$ for $i=0,1,2,3,4$ to $A_k$ and from  $B_k$ as follows,
\bea
S_k (t_1, t_2, t_3 |\tau_2) = \sum _{i=0}^4 S_k^{(i)}  (t_1, t_2, t_3 |\tau_2) 
\eea
We shall  spell out their precise contributions in the subsequent subsections.

\subsubsection{Contributions to $S_k^{(0)}$}

The term $S_k^{(0)} $ arises from the contribution to $A_k$ of $N =0$, and is non-zero only for the constant Fourier mode $k=0$.   Performing a Poisson resummation in $M$ gives,
\bea
S_k ^{(0)} (t_1, t_2, t_3 |\tau_2)  =  \delta _{k,0} \sum _{M \in \ZZ^2}  
 \exp \left \{    - { \pi \over \tau_2}  M^t T M   \right \}
\eea
Partitioning the summation over $M$ according to (\ref{part}), we have,
\bea
S_k ^{(0)} (t_1, t_2, t_3 |\tau_2) 
=   \delta _{k,0} \Big ( 1 + \sum _{i=1}^3 L(t_i' |\tau_2)   + \tilde S_0 ^{(0)} (t_1, t_2, t_3 |\tau_2)    \Big )
\eea
where we have defined $t_i' = t_1+t_2+t_3-t_i$ while the functions $L$ and $\tilde S_0 ^{(0)}$ are defined by,
\bea
\label{S0}
L (t_i'|\tau_2) & = & \sum _{m \in \ZZ} ' e^{- \pi t_i' m^2/\tau_2} 
\no \\
\tilde S_0 ^{(0)} (t_1, t_2, t_3 |\tau_2)  & = &  \sum _{M \in \mN^{(4)}}  
 \exp \left \{    - { \pi \over \tau_2}  M^t T M   \right \}
\eea
The contributions of the functions $L(t)$ naturally combine with those from  $B_k(t)$.

\subsubsection{Contributions to $S_k^{(i)}$ for $i=1,2,3$}

For $i=1,2,3$, the term $S_k^{(i)}$ arises from the contribution to $A_k$ of the partition with $N \in \mN^{(i)}$ and  
 from $B_k(t_i'|\tau_2)$. It will be convenient to exhibit it  as follows,
\bea
S_k^{(i)} (t_1, t_2, t_3 |\tau_2) = \tilde S_k^{(i)} (t_1, t_2, t_3 |\tau_2) 
- \delta _{k,0}  - \delta _{k,0} L(t_i'|\tau_2) 
\eea
where the reduced Fourier mode $\tilde S_k^{(i)} (t_1, t_2, t_3 |\tau_2) $ is given by,
\bea
\label{Si}
\tilde S_k^{(i)} (t_1, t_2, t_3 |\tau_2)  & = & 
 { \tau_2 \over (\det T)^\half } \sum _{N \in \mN^{(i)} } \, \sum _{M \in \ZZ^2 } \!  \delta _{M^t \ep N, k} \, 
 \exp \left \{    - \pi  \tau_2  {M^t T M \over \det T}   - \pi \tau_2  N^t  T  N   \right \} 
\no \\ &&
- B_k (t_i') + \delta _{k,0} \Big ( 1 + L(t_i'|\tau_2) \Big )
 \eea
The  contributions individually fail to exponentially decay in the direction $t_i \to \infty$ but the sum of the two lines above produces a function $S_k^{(i)} $ which exponentially decays in all directions $t_1, t_2, t_3$ at infinity.  To establish uniform exponential decay, it will be convenient to treat the cases $k=0$ and $k \not=0$ separately. 

\sm

For $k=0$, consider the case $i=1$, the other cases being obtained by cyclic permutations of $t_1, t_2, t_3$. For $i=1$, the partition $\mN^{(i)}$ may be parametrized explicitly by $N= (0, n_2)$ with $n_1=0$ and $n_2=-n_3 \not=0$. The constraint $M^t \ep N=0$  reduces to $m_1n_2=0$ so that we must have $m_1=0$, while $m_2 \in \ZZ$. 
Furthermore, the summation over $m,n$ in the  function $B_0(t_2+t_3 |\tau_2)$ is constrained by $mn=0$ and simplifies as follows,
\bea
B_0 (t_2 + t_3 |\tau_2) = \sqrt{ { \tau _2 \over t_2 + t_3} } \sum _{ m \in \ZZ} e^{- \pi \tau_2 m^2/(t_2 + t_3)} + \sqrt{ { \tau _2 \over t_2 + t_3} }
\sum _{n \not=0} e^{- \pi \tau_2 (t_2 + t_3)  n^2}
\eea
Poisson resummation over $m$ in the first term on the right side  gives, 
\bea
B_0 (t_2 + t_3  |\tau_2) = 1 + L(t_2 + t_3 |\tau_2) + \sqrt{ { \tau _2 \over t_2 + t_3 } }
\sum _{n \not=0} e^{- \pi \tau_2 (t_2 + t_3 ) n^2}
\eea
Using the parametrization of $M$ and $N$ for $k=0$ and $i=1$ given above, the sums over $m_2$ and $n_2$ on the first line of the right side of (\ref{Si}) factorize, and the expression takes the form, 
\bea
 { \tau_2 \over (\det T)^\half } \sum _{m_2  } \, \sum _{n _2 \not=0 } \,
 \exp \left \{    - \pi  \tau_2   { t_2 + t_3 \over \det T} \, m_2^2   - \pi \tau_2   (t_2+t_3)  \, n_2^2 \right \} 
 \eea
 Poisson resumming over $m_2$, and combining the expression with the result obtained for $B_0$ shows that the $m_2=0$ mode of the Poisson resummation is cancelled by the terms from $B_0$ and $L$, and gives the following expression for $\tilde S_0^{(1)}$,
 \bea
 \label{Sif}
 \tilde S_0^{(1)} (t_1, t_2, t_3 |\tau_2)  =
\sqrt{ { \tau _2 \over t_2+t_3} }  \sum _{m_2, n_2 \not= 0  } 
 \exp \left \{    - \pi  { \det T \over \tau_2 \,(t_2+t_3)} \, m_2^2   - \pi \tau_2  (t_2+t_3) \, n_2^2 \right \} 
 \eea
 The functions $\tilde S_0^{(2)}$ and $\tilde S_0^{(3)}$ are obtained by cyclic permutations in the variables $t_1, t_2, t_3$ of  (\ref{Sif}), and the resulting functions $\tilde S_0^{(i)} (t_1, t_2, t_3 |\tau_2)$ have uniform exponential decay in all directions of $t_1, t_2, t_3$. The case $k \not=0$ may be handled similarly, but will not be needed to prove the Theorems, and we shall not discuss it further.

\subsubsection{Contributions to $S_k^{(4)}$}

The term $S_k^{(4)} $ arises solely from the contribution to $A_k$ of the partition with $N \in \mN^{(4)}$,  and is given as follows,
\bea
S_k ^{(4)}  (t_1, t_2, t_3 |\tau_2 ) =  { \tau_2 \over (\det T)^\half } \sum _{N \in \mN ^{(4)}} \, \sum _{M \in \ZZ^2 } \!  \delta _{ M^t \ep N, k}  \exp \left \{    - \pi  \tau_2  {M^t T M \over \det T}   - \pi \tau_2  N^t  T  N   \right \}
\eea
Clearly, $S_k^{(4)} $ is uniformly exponentially decaying in $t_1, t_2, t_3$ at infinity. For $k=0$ the constraint $M^t \ep N=0$ forces $M$ to either vanish, or to belong to $\mN^{(4)}$, excluding the cases $M \in \mN^{(i)}$ for $i=1,2,3$. Therefore, it will be convenient to split the sum accordingly,
\bea
S_k ^{(4)}  (t_1, t_2, t_3 |\tau_2 ) = \tilde S_k ^{(4)}  (t_1, t_2, t_3 |\tau_2 ) + \tilde S_k ^{(5)}  (t_1, t_2, t_3 |\tau_2 )
\eea
where
\bea
\label{S45}
\tilde S_k ^{(4)}  (t_1, t_2, t_3 |\tau_2 )
& = & 
{ \tau_2 \delta _{k,0} \over (\det T)^\half } \sum _{N \in \mN ^{(4)}}  \exp \left \{     - \pi \tau_2  N^t  T  N   \right \}
\\
\tilde S_k ^{(5)}  (t_1, t_2, t_3 |\tau_2 )
& = & 
{ \tau_2 \over (\det T)^\half } \sum _{N \in \mN ^{(4)}} \, \sum _{M \not= 0 } \!  \delta _{ M^t \ep N, k}  \exp \left \{    - \pi  \tau_2  {M^t T M \over \det T}   - \pi \tau_2  N^t  T  N   \right \}
\no
\eea
where $\tilde S_k^{(4)}$ arises from $M=0$ while $\tilde S_k^{(5)}$ arises from $M \in \mN^{(4)}$.

\subsection{Summary of contributions to the constant Fourier mode}

Collecting the contributions  obtained in (\ref{S0}), (\ref{Sif}) and (\ref{S45}), we find that the constant Fourier mode is given by the sum of six terms,
\bea
\label{4a1}
C_{a_1, a_2, a_3}^{(0)} ( \tau_2) = \sum _{i=0}^5 \cC _0 ^{(i)} (\tau_2)
\eea
each of which is given by the following integrals over $t_1, t_2, t_3$ of the corresponding functions $\tilde S_0 ^{(i)}$ evaluated in the previous section, 
\bea
\label{4a2}
\cC_0^{(i)} (\tau_2) = \left ( \prod _{r=1}^3 {1 \over \Gamma (a_r)} \int _0 ^\infty dt_r \,  t_r^{a_r-1}  \right ) \tilde S_0^{(i)}(t_1, t_2, t_3 |\tau _2 )
\eea
Each integrand uniformly  decays to zero exponentially fast in any direction as $t_r \to \infty$. When no confusion is expected to arise, we shall often suppress  the dependence on the parameters $a_1, a_2, a_3$ to save notation.

\section{The Laurent polynomial}
\setcounter{equation}{0}
\label{sec:4}

In this section, we shall obtain the Laurent polynomial $\cL(\tau_2)$ in the constant Fourier mode $C_{a_1, a_2, a_3}^{(0)} (\tau_2)$   of the modular functions $C_{a_1, a_2, a_3} (\tau)$, and prove  Theorem \ref{t1}. To this end, we evaluate the contributions $\cC_0^{(i)}(\tau_2)$ for $i=0,1,\cdots, 5$ in the subsections below. The remaining exponential contributions $\cE(\tau_2)$ to the constant Fourier mode will be evaluated in the subsequent section.

\subsection{Evaluating $\cC_0 ^{(0)}$} 

The integral over $S_0^{(0)} $ evaluates to a sum over $M \in \mN^{(4)}$ which may be parametrized by,
\bea
\cC_0 ^{(0)}  (\tau_2) = { \tau_2^w \over \pi^w}   \sum _{m_1 \not=0} { 1 \over m_1 ^{2a_1} } \sum _{m_2 \not= 0, -m_1} { 1 \over m_2^{2 a_2} (m_1+m_2)^{2a_3}}
\eea
To compute the infinite sum over $m_2$, we proceed by decomposing the summand into partial fractions in $m_2$, using the general partial fraction decomposition formulas,  valid for $a, b \in \NN$,
\bea
\label{pf1}
{ 1 \over (z+x)^a (z+y)^b} = \sum _{k=1}^a { \cA _k(a,b) \over (z+x) ^k\, (y-x)^{a+b-k} } + \sum _{k=1}^b { \cB _k(a,b)  \over (z+y)^k \, (y-x)^{a+b-k} }
\eea
where $\cA_k(a,b)$ and $\cB_k(a,b)$ are given by binomial coefficients,
\bea
\label{pf2}
\cA_k(a,b) & = &  (-)^{a+k} \binom{a+b-k-1}{a-k}
\no \\
\cB_k(a,b) & = &  (-)^a \binom{a+b-k-1}{b-k}
\eea
For the case at hand, we set $x=0$, $y=m_1$,  $z=m_2$, for the positive integers exponents $a=2a_2$ and $b=2a_3$, and we find, 
 \bea
\sum _{m_2 \not= 0, -m_1} 
{ 1 \over m_2^{2a_2} (m_1+m_2)^{2a_3}}
& = &
 \sum _{k=1} ^{ a_2} { 2 \zeta (2k)  \over m_1^{2a_2+2a_3-2k} } \Big  ( \cA_{2k}(2a_2,2a_3) + 
 \cB_{2k}(2a_2,2a_3) \Big )
\no \\ && 
- { 1 \over m_1^{2a_2+2a_3} }{ \Gamma (2a_2+2a_3+1) \over \Gamma (2a_2+1) \Gamma (2a_3+1) }
\eea 
Evaluating next the sum over $m_1$, and expressing the resulting even $\zeta$-values in terms of Bernoulli numbers using,
\bea
\zeta (2k) =  \half (2 \pi)^{2k} (-)^{k+1} { B_{2k} \over (2k)!}
\eea
we find,  
\bea
\label{4b4}
 \cC_0^{(0)} (\tau_2)
= 
 (-4 \pi \tau_2)^w \sum _{k=0} ^{ a_2}   {  B_{2k} B_{2w-2k} \over (2k)! \, (2w-2k)! } \, { \Gamma (2a_2+2a_3-2k) \over \Gamma (2a_3) \Gamma (2a_2-2k+1)}
+ (a_2 \leftrightarrow a_3)
\eea
This formula reproduces correctly the top terms previously evaluated in (\ref{spec}).

\subsection{Evaluating $\cC_0 ^{(i)}$ for $i=1,2,3$} 

We  evaluate the case $i=3$, the cases $i=1,2$ being obtained by cyclic permutations of $a_1, a_2, a_3$. The integral over $t_3$ in $S_0^{(3)}$ may be readily carried since the dependence of $S_0^{(3)}$ on $t_3$ is entirely contained in $\det T$, whose dependence on $t_3$ is as follows,
\bea
{ \det T \over t_1+t_2} = t_3 + { t_1 t_2 \over t_1+t_2}
\eea
Carrying out the integral over $t_3$ gives, 
\bea
\cC_0^{(3)} (\tau_2) & = &
  \sum _{m, n \not=0} \left ( { \tau_2 \over \pi m^2} \right )^{a_3}
\prod _{r=1}^2  \int _0 ^\infty {dt_r \, t_r^{a_r-1} \over \Gamma (a_r)}  \sqrt{{ \tau_2 \over t_1 + t_2}}\,
\no \\ && \hskip 1in \times
\exp \left \{ - \pi \tau _2 n^2 (t_1+t_2)  - \pi {m^2 \over \tau_2} \,  {t_1 t_2  \over t_1+t_2} \right \}
 \eea
Parametrizing the integration variables by $t_1= x t$ and $t_2 = t (1-x) $ with $t\geq 0$ and $  0 \leq x \leq 1$, and carrying out the  integral in $t$ produces the following result, 
\bea
\label{C3G}
\cC_0^{(3)} (\tau_2) = { \tau_2^w \over \pi^w}
\sum _{m,n \not = 0} { 1 \over |m|^{2w-1}} \, G_{a_1, a_2} \left ( {  \tau_2 n \over m } \right )
\eea
where the function $G_{a_1,a_2}(\mu)$ is given by the integral representation, 
\bea
\label{G1}
G_{a_1,a_2}(\mu) = {\sqrt{\pi}  \, \Gamma (a_1+a_2-\half)  \over  \Gamma (a_1) \Gamma (a_2) }
\int _0 ^1 dx \,  { x ^{a_1-1} (1-x)^{a_2-1}  \over  \big ( \mu^2 +  x(1-x) \big ) ^{a_1+a_2-\half}}
 \eea
The function $G_{a_1,a_2}(\mu)$ is even in $\mu$, and invariant under interchanging $a_1$ and $a_2$. To evaluate it, we shall make use of the following Lemma.

\medskip
\noindent
\begin{lem}\label{l1}  The function $G_{a_1, a_2}(\mu)$ defined in (\ref{G1}) admits the equivalent representations.
{\sl \begin{description}
\item (a)  For $\mu \in \RR$, and $a_1, a_2 \in \CC$ with $\Re(a_1), \Re(a_2) \geq 1$, 
\bea
\label{G0}
G_{a_1, a_2} (\mu) = \int _\RR du \, { 1 \over (u^2+\mu^2)^{a_1} ((u+1)^2+\mu^2)^{a_2}}
\eea 
\item (b) For $a_1, a_2 \in \NN$,  
\bea
\label{G}
 G_{a_1, a_2} (\mu)  = 
\left  ( \sum _{\alpha =0}^{a_1-1} \sum _{\beta =0} ^{a_1-1-\alpha} 
{-  i \pi \, g_{a_1,a_2}(\alpha, \beta)  \over  (2i\mu) ^{2a_1-1-\alpha-\beta} (1+2i\mu )^{a_2+\alpha}} + \hbox{c.c.} \right ) + ( a_1 \leftrightarrow a_2) 
\eea
The coefficients $g_{a_1, a_2} (\alpha, \beta)$ are given by the product of binomial coefficients of (\ref{G3}).
\end{description}}
\end{lem}

\sm

To prove part {\sl (a)} of  Lemma 1, we start from expression (\ref{G0}) and use standard techniques for the evaluation of Feynman diagrams in quantum field theory to derive its expression given in (\ref{G1}). One makes use of an integral representation formula for a product of denominators,
\bea
{ 1 \over A_1^{a_1} A_2^{a_2}} 
= { \Gamma (a_1+a_2) \over \Gamma (a_1) \Gamma (a_2)} 
\int _0 ^1 dx { x^{a_1-1} (1-x)^{a_2-1} \over ( x A _1+ (1-x) A_2)^{a_1+a_2}}
\eea
valid for $A_1, A_2 >0$ and $\Re(a_1), \Re(a_2) >0$, and the evaluation of the resulting $u$-integral, 
\bea
\int _\RR du { 1 \over (u^2 +\mu^2)^a} = { \sqrt{\pi} \, \Gamma (a-\half) \over \Gamma (a) \, |\mu|^{2a-1}}
\eea
valid for $\Re(a) > \half$. To prove part {\sl (b)} of Lemma 1, we make use of the fact that, for $a_1, a_2 \in \NN$,  the integrand in (\ref{G0}) is a rational function of $u$, with poles at $u = \pm i \mu$ and $u=-1 \pm i \mu$. The integral may then be evaluated by standard residue methods.

\subsubsection{Calculating the infinite sums for $\cC_0^{(3)}$}

Formula (\ref{C3G}) expresses $\cC_0^{(3)}$ as a double sum over the function $G_{a_1, a_2}$. We use the expression  for $G_{a_1, a_2}$ obtained in (\ref{G}), interchange the order of the  sums over $m,n$ with the sums over $\alpha, \beta$, and combine the factors of $i$ with factors of the absolute values $|m|$, 
\bea
\cC_0^{(3)} = { \tau_2^w \over \pi^w}
\left  ( \sum _{\alpha =0}^{a_1-1} \sum _{\beta =0} ^{a_1-1-\alpha} \sum _{m,n \not = 0}
{(-)^w (-i|m|)^{1-2w} \pi \, g_{a_1,a_2}(\alpha, \beta)  \over  (i|{ 2 \tau_2n \over m} |) ^{2a_1-1-\alpha-\beta} (1+i|{ 2 \tau_2 n \over m} | )^{a_2+\alpha}} + \hbox{c.c.} \right ) 
+ ( a_1 \leftrightarrow a_2) 
\eea
The role of the complex conjugate contribution is to reverse the sign of $-i|m|$, so that we may omit the absolute value symbol on $m$ and include a factor of 2 to account for the addition of the complex conjugate term. Restricting the sum over $n$ to $n >0$ gives another factor of 2, and rearranging the factors of $i$, we express $\cC_0^{(3)}$ as follows,
\bea
\label{CKK}
\cC_0^{(3)}  =   \cK(a_1,a_2,a_3) + \cK(a_2,a_1,a_3)
\eea
where $\cK$ is given by,
\bea
\cK(a_1,a_2,a_3) = 
{ \tau_2^w \over \pi^w}
 \sum _{\alpha =0}^{a_1-1} \, \sum _{\beta =0} ^{a_1-1-\alpha} \,
\sum _{n=1}^\infty \,  \sum _{m \not = 0} \, { - 4 \pi i  \, g_{a_1,a_2}(\alpha, \beta) 
 \over m^{A} (m+2 i \tau_2 n  )^{B} ( 2 i \tau_2n )^{C} }
\eea 
where we have used the following abbreviations, 
\bea
\label{ABC}
a& = & a_2 + 2 a_3 +\beta
\no \\
b& = & a_2 + \alpha
\no \\
c& = & 2 a_1 -\alpha - \beta -1
\eea
The sum over $m$ may be carried out using the partial fraction decomposition formulas of (\ref{pf1}) and (\ref{pf2}) for the parameters $x=0$, $y=2i \tau _2 n$ with $n>0$, $z=m$, and the exponents $a$, $b$ defined in (\ref{ABC}), and we find, 
\bea
\sum _{ m \not= 0} { 1 \over m^a(m+2i \tau_2 n)^b} 
& = & \sum _{k=1}^{[a/2]}  { 2 \zeta (2k) \, \cA _{2k} (a,b)  \over (2i \tau_2 n)^{a+b-2k} }   
- \sum _{k=1}^b { \cB _k (a,b) \over  (2i \tau_2 n)^{a+b} }
- { i \pi \cB_1 (a,b) \over (2i \tau_2 n)^{a+b-1}}
\no \\ &&
+ \sum _{k=1}^b { \cB_k(a,b) \over (2i \tau_2 n) ^{a+b-k} }  {(- 2 \pi i)^k \over \Gamma (k) } 
\sum _{p=1}^\infty p^{k-1} e^{-4 \pi  p n \tau_2}
\eea
The sum over $n$ is carried out by multiplying the above relation by  $- 4i (2i \tau_2 n)^{-c}$ with $c$ given in (\ref{ABC}),  and separating the resulting sum into the contributions to the Laurent polynomial part and the contributions  to the exponential part, 
\bea
\cK(a_1,a_2,a_3) =\cK_L(a_1,a_2,a_3) + \cK_E(a_1,a_2,a_3)
\eea
 The Laurent polynomial part is given by, 
\bea
\label{KL}
\cK_L(a_1,a_2,a_3) & = &  
\sum _{\alpha =0}^{a_1-1} \sum _{\beta =0} ^{a_1-1-\alpha} 
g_{a_1,a_2}(\alpha, \beta)  \Bigg [
 { \zeta (2w-2)  \cB_1 (a,b) \over   (-4\pi \tau_2 )^{w-2}  } 
 \no \\ && \hskip 0.5in
 + {  2  \zeta (2w-1) \over  (-4 \pi \tau_2 )^{w-1} } { (-)^{a_2+\beta} \, \Gamma (2a_2+2a_3+\alpha+\beta) \over \Gamma (a_2+2a_3+\beta+1) \Gamma (a_2+\alpha)}
\no \\ && \hskip 0.3in 
- \sum _{k=1}^{a_3+[(a_2+\beta)/2]}  { 4 \, (-)^{k}  \zeta (2k) \, \zeta(2w-1-2k) \cA _{2k} (a,b)  \over (-4 \pi \tau_2 )^{w-1-2k} \, (2 \pi)^{2k} } 
 \Bigg ]  
\eea 
while the purely exponential part is given by,  
\bea
\cK_E(a_1,a_2,a_3)  =   
 \sum _{\alpha =0}^{a_1-1} \sum _{\beta =0} ^{a_1-1-\alpha} \!
g_{a_1,a_2}(\alpha, \beta)  \sum _{k=1}^{a_2+\alpha} 
{ 2 \, (-)^w \, \cB_k (a,b) \over (4 \pi \tau_2 ) ^{w-1-k} \,  \Gamma (k)}   
\sum _{n=1}^\infty n^{k-1} \sigma _{2-2w} (n) q^n \bar q^n 
\eea
To simplify the purely exponential part, we have used the  standard rearrangement formula, 
\bea
\sum _{n=1}^\infty { 1 \over n^{2w-1-k}} \sum _{p=1}^\infty p^{k-1} e^{-4 \pi  p n \tau_2}
=
\sum _{n=1}^\infty n^{k-1} \sigma _{2-2w} (n) q^n \bar q^n
\eea

\subsection{Evaluating $\cC_0 ^{(4)}$} 

We shall use the following integral representation for the  factor $\det T$ in $\tilde S_0 ^{(4)}$,
\bea
{ 1 \over \tau_2 (\det T)^\half}
= \int _{\RR^2} du dv \, e^{  - \pi \tau_2 \left ( t_1u^2 + t_2 v^2 + t_3 (u+v)^2 \right ) }
\eea
in order to decouple the $t$-integrals in (\ref{4a2}) for this function, and we obtain, 
\bea
\label{4c2}
\cC_0^{(4)} =  { \tau_2^2 \over (\pi \tau_2)^w} \sum _{N \in \, \mN^{)4)} } \, \int _{\RR^2} 
 { du dv \over (u^2+n_1^2)^{a_1} (v^2+n_2^2)^{a_2} ((u+v)^2+n_3^2)^{a_3} }
\eea
The integral is independent of  $\tau_2$ so that $\cC_0 ^{(4)}$ contributes exclusively to the order $\tau_2^{2-w}$  in the Laurent polynomial. To evaluate the integral over $\RR^2$ and the summation over $\mN^{(4)}$, we shall proceed as follows. 

\sm

We begin by simplifying the summation over $\mN^{(4)}$. Invariance of the set $\mN^{(4)}$ under permutations of $n_1,n_2,n_3$ guarantees invariance of $\cC_0 ^{(4)}$ under permutations of $a_1, a_2, a_3$. We partition $\mN^{(4)}$ into three disjoint subsets,  $n_1n_2>0$, $n_2 n_3>0$ and $n_1n_3>0$, and we may restrict the summation over $\mN^{(4)}$ to any single one of these subsets provided we add the contribution of the two cyclic permutations of $a_1,a_2, a_3$. We shall choose the subset $n_1 n_3>0$. Within this subset, the sectors $n_1, n_3>0$ and $n_1, n_3 <0$ contribute equally; we shall restrict to $n_1, n_3 \geq 1$ upon including a factor of 2. 

\sm

Since the exponents $a_1, a_2, a_3$ are positive integers, the integrals  are over a rational function $f(u,v)$, and may be evaluated by residue methods. Taking the above preparations into account, we obtain the following expression, 
\bea
\cC_0^{(4)} & = &  { 2 \tau_2^2 \over (\pi \tau_2)^w} \sum _{n_1, n_3 =1}^\infty  \, \int _{\RR^2} du \, dv \,  | f (u,v)|^2 
+ \hbox{2 cyclic permutations of } a_1, a_2, a_3
\no \\
f(u,v) & = & { 1 \over (u+in_1)^{a_1} (v+in_2)^{a_2} (u+v+in_1+in_2)^{a_3} }
\eea
The function $f$ is not unique and is chosen such that the last denominator argument is the sum of the preceding two. This choice guarantees that $f(u,v)$ will have a convenient partial fraction expansion in the variable $u$, the integral over which we shall carry out first, 
\bea
f(u,v)  =  
\sum _{k_1=1} ^{a_1} { (-)^{a_1-k_1}  \binom{a_1+a_3-k_1-1}{a_3-1}
\over (u+in_1)^{k_1} (v+in_2)^{w-k_1} }
+ \sum _{k_3=1} ^{a_3}  { (-)^{a_1}   \binom{a_1+a_3-k_3-1}{a_1-1} \over (u+v-in_3)^{k_3} (v+in_2)^{w-k_3} }
\eea
Because the sum is restricted to $n_1 , n_3 >0$, the $u$-integrals in the cross terms of the product $f(u,v) \, \overline{f(u,v)}$ vanish identically since in each case the poles are all either in the upper half plane or all in the lower half plane. Hence  only the integrals in the direct terms contribute, 
\bea
\mI_{1} (n_1, n_3) & = & { 2^{2w} \over 16 \pi^2} 
\int _{\RR^2}  { (-)^{k_1+\ell_1} \, du \, dv \over (u+in_1)^{k_1}(u-in_1)^{\ell_1} (v+in_2)^{w-k_1} (v-in_2)^{w-\ell_1}}
\no \\
\mI_{2} (n_1, n_3) & = & { 2^{2w} \over 16 \pi^2} 
\int _{\RR^2}  { du \, dv \over (u+v-in_3)^{k_3}(u+v+in_3)^{\ell_3} (v+in_2)^{w-k_3} (v-in_2)^{w-\ell_3}}
\eea
In terms of these integrals, $\cC_0^{(4)}$ is given by,
\bea
\cC_0^{(4)} & = &  
~ {2 \over (4 \pi \tau_2)^{w-2} }  \sum_{k_1=1}^{a_1} \sum _{\ell_1=1}^{a_1} 
\binom{a_1+a_3-k_1-1}{a_3-1}  \binom{a_1+a_3-\ell_1-1}{a_3-1} \sum _{n_1,n_3=1}^\infty \mI_{1} (n_1, n_3)
\no \\ && +
 {2 \over (4 \pi \tau_2)^{w-2} } \sum_{k_3=1}^{a_3} \sum _{\ell_3=1}^{a_3} 
  \binom{a_1+a_3-k_3-1}{a_1-1}  \binom{a_1+a_3-\ell_3-1}{a_1-1} \sum _{n_1,n_3=1}^\infty  \mI_{2} (n_1, n_3)
 \no \\ &&
 + \hbox{ 2 cyclic permutations of } a_1, a_2, a_3
\eea
The integrals over $u$ and $v$ are manifestly decoupled from one another in $\mI_1$, as well as  in $\mI_{2}$  after performing the shift $u \to u-v$. Using residue methods, the integrals evaluate to,
\bea
\mI_{1} (n_1, n_3)& = & 
\binom{k_1+\ell_1-2}{k_1-1} \binom{2w-k_1-\ell_1-2}{w-k_1-1}
{ \theta(-n_1n_2)  + (-)^{k_1+\ell_1}  \theta(n_1n_2) \over \, |n_1|^{k_1+\ell_1-1} \, \, |n_2|^{2w-k_1-\ell_1-1}}
\no \\ 
\mI_{2} (n_1, n_3)  & = & 
\binom{k_3+\ell_3-2}{k_3-1} \binom{2w-k_3-\ell_3-2}{w-k_3-1} 
{ \theta(-n_2n_3)  + (-)^{k_3+\ell_3}  \theta(n_2 n_3) \over \, |n_3|^{k_3+\ell_3-1} \, \, |n_2|^{2w-k_3-\ell_3-1}}
\eea
The contributions of the integrals $\mI_{1}$ and $\mI_{2}$ to $\cC_0^{(4)}$ are related by permuting $a_1$ and $a_3$. Hence we may retain only the sum involving $\mI_{1}$ provided we then include all five permutations of $a_1, a_2, a_3$. Putting all together we have,
\bea
\label{C4Z}
\cC_0^{(4)}  =   
~ {2 \over (4 \pi \tau_2)^{w-2} }  \sum _{\sigma \in \mS_3} Z(a_{\sigma (1)} ,a_{\sigma (2)} ,a_{\sigma (3)} ) 
\eea
where the function $Z(a_1,a_2,a_3)$ was defined in (\ref{ZZ}).

\subsection{Evaluating $\cC_0^{(5)}$}

Combining the second equation of (\ref{S45}) with (\ref{4a2}) we see that $\cC_0^{(5)}$ is given by, 
\bea
\label{5a1}
\cC_0^{(5)}  = 
 \prod _{r=1}^3  \int _0 ^\infty dt_r \,{ t_r^{a_r-1} \over \Gamma (a_r)}  
 \sum _{M, N \in \, \mN^{(4)}}   { \tau_2 \, \delta _{ M^t \ep N, 0}  \over (\det T)^\half } 
\, \exp \left \{ - \pi  \tau_2   {M^t T M \over \det T}   - \pi \tau_2  N^t  T  N   \right \}
\eea
\sm
We begin by parametrizing the space of matrices $M,N$. To  satisfy the condition $M^t \ep N=0$, the column matrices $M,N$ must be proportional to one another, and thus proportional to a common matrix $K$ with integer entries $k_1, k_2$. Since $M,N \in \mN^{(4)}$, the numbers $m_i, n_i, k_i$ with $i=1,2,3$ 
are also non-vanishing.  The complete  solution is given by,
\bea
\label{5a2}
M = \mu K \hskip 0.6in N= \nu K \hskip 0.6in K = \binom{k_1} {k_2} 
\eea
where we choose $k_1$ and $ k_2$ relatively  prime, $k_1 >0$, and  $\mu, \nu \in \ZZ$ with $\mu, \nu \not=0$.   Since the summand in $\cC_0^{(5)}$ depends only on $\mu^2$, $\nu^2$, $k_1^2, k_2^2$ and $ k_3^2$, we may include in the sum over $k_1$ and $k_2$ pairs with arbitrary signs along with an overall factor of $\half$, and restrict $\mu, \nu$ to be positive upon including an overall factor of 4. Finally, we may choose a particular ordering of $k_1^2, k_2^2 , k_3^2$ upon including  symmetrization under all six permutations of the variables $a_1, a_2, a_3$. We shall denote the space of such pairs $k_1, k_2$ of relatively prime, ordered, integers by $\mK$.

\sm

Next, we change integration variables from $(t_1, t_2, t_3)$ to $(x_1, x_2, t_3)$ with $t_1=x_1 t_3$ and $t_2=x_2t_3$ for $x_1, x_2 \geq 0$, and subsequently change variables from $(x_1, x_2, t_3)$ to $(x_1, x_2, t)$ with,
\bea
t_3 = \left | { \mu \over \nu} \right | { t \over \sqrt{x_1+x_2+x_1x_2}}
\eea
In terms of these new variables, the integrals become,
\bea
\label{5a3}
\cC_0^{(5)} & = &
2 \tau_2 \sum _{(k_1, k_2) \in \mK} \, \sum _{\mu, \nu =1}^\infty  \left ( {\mu \over \nu} \right )^{w-1} 
\int _0^\infty dt \, {t^{w-2} \over \Gamma (a_3)}  \prod _{r=1}^2  \int _0 ^\infty dx_r \,{ x_r^{a_r-1} \over \Gamma (a_r)}  { 1\over (x_1+x_2+x_1x_2)^{w/2}}
\no \\ &&  \times 
 \exp \left \{ - \pi \tau_2 \mu \nu \left (t + { 1 \over t} \right ) { x_1 k_1^2 + x_2 k_2^2 + k_3^2 \over \sqrt{x_1+x_2+x_1x_2}} \right \} + \hbox{ 5 perms of } a_1, a_2, a_3
\eea
The Mellin transform of $\cC_0 ^{(5)}$  with respect to $\tau_2 $ is defined by,
\bea
\cM \cC_0 ^{(5)} (s) = \int _0 ^\infty d \tau _2 \, \tau _2 ^{s-1} \cC_0 ^{(5)}(\tau_2)
\eea
Carrying out the integration over $t$, and performing the sums over $\mu$ and $\nu$ we obtain, 
\bea
\cM \cC_0 ^{(5)} (s)  =  
{ \xi (s+w) \xi (s-w+2) \over  \Gamma (a_1) \Gamma (a_2) \Gamma (a_3)} \sum _{(k_1, k_2) \in \mK} \, 
\prod _{r=1}^2 \int _0^\infty dx_r \, x_r ^{a_r-1} { (x_1+x_2+x_1x_2)^{{s-w+2\over2}} 
\over (x_1 k_1^2 + x_2 k_2^2 + k_3^2)^{s+1}}
\eea
The function  $\xi(s) = \Gamma (s/2) \zeta (s)$ has single poles at $s=0$ and $s=1$ and is analytic elsewhere.
We now prove the following Lemma.

\sm
\noindent
\begin{lem}\label{l4}
Near the cusp $\tau_2 \to  \infty$, the function $\cC_0^{(5)}$ decays exponentially, 
\bea
\cC_0 ^{(5)} (\tau_2) = \cO(e^{-4 \pi \tau_2})
\eea
up to factors which are power-behaved in $\tau_2$, controlled by the poles in the Mellin transform.
\end{lem}
To prove the Lemma, we shall extract out of the integrand a factor of $e^{-4 \pi \tau_2}$. First, we have the standard sharp bound $t + t^{-1} \geq 2$. Secondly, we shall show that the last factor of the exponential has a uniform  sharp lower bound by $2k_1k_2$. To prove it, we set $y_1=x_1+1$ and $y_2=x_2+1$, with $y_1, y_2 \geq 1$,  in terms of which the expression becomes, 
\bea
{ x_1 k_1^2 + x_2 k_2^2 + k_3^2 \over \sqrt{x_1+x_2+x_1x_2}}
= { y_1 k_1^2 + y_2 k_2^2 + 2 k_1k_2 \over \sqrt{y_1y_2-1}}
\eea
Since $(k_1, k_2 ) \in \mK$ we have $k_1, k_2>0$.  Keeping the product $y_1y_2$ fixed, the combination $ y_1 k_1^2 + y_2 k_2^2$ is bounded from below by $2 k_1k_2 \sqrt{y_1y_2}$. The resulting lower bound equals $2k_1k_2$ times a function of $y_1y_2$ which is independent of $k_1, k_2$ and whose lower bound is 1. Thus, for  $(k_1, k_2) \in \mK$, we may extract out of the integrand a factor of $e^{-4 \pi \tau_2 \mu \nu k_1 k_2}$. The Lemma follows from  $\mu, \nu , k_1k_2\geq 1$, and the observation that the remaining integrations converge.
Their power behavior is governed by the positions of the poles of the Mellin transform.

\subsection{Proof of Theorem \ref{t1}}

Having shown that $\cC_0^{(5)}$ does not contribute to the Laurent polynomial $\cL$, Theorem \ref{t1} may now be proven by collecting the contributions to $\cL$ from the calculations of $\cC_0^{(i)}$ for $i=0,1,2,3,4$ performed in the preceding subsections. From the explicit results in (\ref{4b4}),  (\ref{CKK}), (\ref{KL}) and (\ref{C4Z}), we deduce the following results. 
\begin{enumerate}
\item 
$c_w$ receives contributions exclusively from $\cC_0^{(0)}$ which was evaluated explicitly in (\ref{4b4}). The coefficient $c_w$ given by (\ref{4b4}) is manifestly a rational number. This proves  part (a). 
\item 
$c_{w-2k-1}$ for $1\leq k \leq w-1$ receives contributions exclusively from $\cC_0^{(i)}$ with $i=1,2,3$, obtained from  the second and third terms under the sum in (\ref{KL}). The coefficients $c_{w-2k-1}$ given by (\ref{4b4}) are manifestly  rational numbers. This proves part~(b). 
\item 
$c_{2-w}$ receives contributions exclusively from the first term under the sum in (\ref{KL}) which gives rise to the term in $\zeta(2w-2)$ in (\ref{C2W}), as well as from the entire contribution of $\cC_0^{(4)}$ in (\ref{C4Z}).
Putting both together proves part (c).
\end{enumerate}

\section{Differential equations and exponential terms}
\setcounter{equation}{0}
\label{sec:5}

For given weight $w=a_1+a_2+a_3$ the functions $C_{a_1,a_2,a_3}(\tau)$ with $a_1, a_2, a_3 \in \NN$   satisfy a system of inhomogeneous linear differential equations whose inhomogeneous part is a linear combination of the non-holomorphic Eisenstein series $E_w$ and products of the form $E_{w-\ell} E_\ell$ with $ 2 \leq \ell \leq w-2$. Analyzing the general structure of these differential equations and their solutions provides convenient paths towards proving Theorems \ref{t2} and \ref{t3}, which we shall carry out in the present section.

\subsection{Inhomogeneous Laplace-eigenvalue equations}

The two-loop modular graph functions satisfy a system of differential equations, 
\bea
\label{diff}
\left (  \Delta -  \sum _{i=1}^3 a_i(a_i-1)   \right ) C_{a_1, a_2, a_3} & = &  
a_1a_2 \Big ( C_{a_1-1,a_2+1,a_3} + \half C_{a_1+1,a_2+1,a_3-2} - 2 C_{a_1,a_2+1,a_3-1} \Big ) 
\no \\ && + \hbox{ 5 permutations of } a_1, a_2, a_3 
\eea
valid for $a_r\geq 3$ for $r=1,2,3$. The validity of these equations may be extended to $a_r \geq 1$ by supplementing them with the following degenerate cases,
\bea
\label{supp}
C_{a_1,a_2,0} & = & E_{a_1} \, E_{a_2} - E_{a_1+a_2} \hskip 1.5in a_1+a_2 \geq 3
\no \\
C_{a_1,a_2,-1} & = & E_{a_1-1} E_{a_2} + E_{a_1} E_{a_2-1} \hskip 1.2in a_1,a_2 \geq 2
\eea
The right side of (\ref{supp}) may involve the symbol $E_1$ (formally corresponding to a divergent series), but its contribution systematically cancel out of the right side of (\ref{diff}). For example, the lowest weight cases are as follows,
\bea
\Delta \, C_{1,1,1} & = & 6 \, E_3
\no \\
\Delta \, C_{2,1,1} & = & 2 \, C_{2,1,1} + 9 \, E_4 - E_2^2
\no \\
\Delta \, C_{2,2,1} & = & 8 \, E_5
\eea
The Laplacian preserves the weight $w=a+b+c$ of the modular graph functions, and the equations of (\ref{diff}) and (\ref{supp}) may be viewed as acting on the space of modular graph functions of given weight $w$, provided we  assign the weight $s$ to the Eisenstein series $E_s$.

\sm

The  Fourier series for a non-holomorphic Eisenstein series $E_n(\tau)$ with integer $n \geq 2$, was given in (\ref{2d5}). It was shown in \cite{D'Hoker:2015foa} that the Laplace operator $\Delta = 4 \tau_2^2 \p_\tau \p_{\bar \tau}$  on the space of functions $C_{a_1,a_2,a_3}(\tau)$ for weight $w=a_1+a_2+a_3 \geq 3$ may be ``diagonalized" resulting in ``eigenfunctions" $\cC_{w,s,\mpp} (\tau)$ which are linear combinations of the functions $C_{a_1, a_2,a_3}(\tau)$ of weight $w$ which obey the following type of equation, 
\bea
(\Delta - s(s-1) ) \cC_{w,s,\mpp} (\tau) = \cH_0(w,s,\mpp) E_w(\tau) 
+ \sum _{\ell=2}^{w-2} \cH_\ell(w,s,\mpp)  E_{w-\ell}(\tau) E_\ell(\tau)
\eea
where $s$ and $\mpp$ are integers running over the following ranges,
\bea
\label{smp}
s= w-2\mm  
\hskip 0.8in  
1 \leq \mm \leq \left [ { w-1 \over 2} \right ]
\hskip 0.8in 
0 \leq \mpp \leq \left [ { s-1 \over 3} \right ]
\eea
The coefficients $\cH_\ell$ are combinatorial rational numbers which depend on $a,b,c,s,\mpp,\ell$. Decomposing $\cC_{w,s,\mpp}(\tau)$ into Fourier modes, 
\bea
\cC_{w,s,p}(\tau) = \sum _{k \in \ZZ} \cC_{w,s,\mpp}^{(k)} (\tau_2) \, e^{ 2 \pi i k \tau_1}
\eea
each Fourier mode satisfies a separate ODE in $\tau_2$,
\bea
&&
\Big (\tau_2^2 \p_{\tau_2}^2  - 4 \pi^2 k^2 \tau_2^2 - s(s-1) \Big ) \cC_{w,s,\mpp}^{(k)} (\tau_2) 
\\  && \qquad =
\cH_0(w,s,\mpp)  \int _0^1 d \tau _1 e^{-2\pi i k \tau_1} E_w(\tau) 
+ \sum _{\ell=2}^{w-2} \cH_\ell (w,s,\mpp)  \int _0^1 d \tau _1 e^{-2\pi i k \tau_1} E_{w-\ell}(\tau) E_\ell(\tau)
\no
\eea 

\subsection{Proof of Theorem \ref{t2}}

To prove Theorem \ref{t2}, we concentrate on the constant Fourier mode for $k=0$ and in particular on its contribution proportional to the functional behavior $c_{2-w} (4 \pi \tau_2) ^{2-w}$. The Fourier transform on $E_w$ makes  vanishing contribution to this power of $\tau_2$, and the contributions of the bilinears $E_{w-\ell} E_\ell$ are readily evaluated using (\ref{2d5}), and we find, 
\bea
&&
(w-s-1)(w+s-2) c_{2-w} 
\no \\ && \qquad
 =
 \sum _{\ell=2}^{w-2} 16 \cH_\ell (w,s,\mpp)  \binom{2w-2\ell-3}{w-\ell-1}  \binom{2\ell-3}{\ell-1} 
 \zeta (2w-2\ell-1)\zeta (2\ell-1)
\eea 
Given that the allowed values of $s$ in the spectral decomposition of the functions $C_{a_1, a_2, a_3}$ onto the basis of eigenfunctions $\mC_{w;s;\mpp}$ is given by (\ref{smp}), and that this decomposition is with rational coefficients, it is clear that the kernel of the operator on the left side of the above equation must vanish. Therefore, $c_{2-w}$ must be a linear combination, with rational coefficients,  of products of pairs of odd $\zeta$-values whose weights sum to $2w-2$., thus proving Theorem \ref{t2}. In particular, any contributions valued in $\pi ^{2w-2}\QQ$ must vanish.

\subsection{Proof of Theorem \ref{t3}}

We shall now concentrate on the constant Fourier mode $\cC_{w,s,\mpp}^{(0)}(\tau_2)$. Having already determined the Laurent polynomial contribution to $C_{a_1, a_2, a_3}(\tau)$ at the cusp we shall consider here its exponential contribution, which we denote by $\cC^{(e)} _{w,s,\mpp}(\tau_2)$. The  contribution to $\cC^{(e)} _{w,s,\mpp}(\tau_2)$ from the function $E_w$ vanishes, so we are left with the contributions from the products $E_{w-\ell} E_\ell$ with $2 \leq \ell \leq w-2$, which obeys the following differential equation,
\bea
\Big (\tau_2^2 \p_{\tau_2}^2  - s(s-1) \Big ) \cC^{(e)} _{w,s,\mpp}(\tau_2)   =  
 \sum _{\ell=2}^{w-2} \cH_\ell(w,s,\mpp) \int _0^1 d \tau _1  E_{w-\ell}(\tau) E_\ell(\tau) \Bigg |_{\hbox{exp}}
\eea 
The right side is readily evaluated and we find, 
\bea
\Big (\tau_2^2 \p_{\tau_2}^2  - s(s-1) \Big ) \cC^{(e)} _{w,s,\mpp}(\tau_2)
& = &  
 \sum _{\ell=2}^{w-2} \cH_\ell (w,s,\mpp) \sum_{n=1}^\infty {  \sigma _{1-2w+2\ell}(n) \sigma _{1-2\ell}(n)  \over (w-\ell+1)! (\ell-1)!} 
 \no \\ && \quad \times
\, 8 \, n ^{w-2} \,  e^{-4 \pi \tau_2 n} P_{w-\ell}(4 \pi \tau_2 n) P_\ell (4 \pi \tau_2 n)
\eea
where the polynomials $P_n$ are given by (\ref{Pw}). 
The general structure of this equation, in terms of the variable $y=4 \pi \tau_2$ is as follows,
\bea
\Big (  y^2 \p_y^2 -s(s-1) \Big ) f_s(y) = \sum _{n=1}^\infty \sum _{m=0 } ^{w-2} f_{m,n} { e^{-ny} \over (ny)^m}
\eea
To solve this equation, we solve for each power of $m$ on the right side of the above equation,
\bea
\Big (  y^2 \p_y^2 -s(s-1) \Big ) \f_{s,m}(y) =  { e^{-y} \over y^m}
\eea
so that $f_s$ is given by, 
\bea
f_s(y) = \sum _{n=1}^\infty \sum _{m=0 } ^{w-2} f_{m,n} \f_{s,m}(ny)
\eea
The solution for $\f_{s,m}(y)$ is given in terms of the incomplete $\Gamma$-function, defined by,
\bea
\Gamma (a,x) = \int _x ^\infty dt \, t^{a-1} \, e^{-t}
\eea
and we find, 
\bea
\f_{s,m}(y) =  { y^s \over 2s-1} \Gamma (-s-m,y) 
+ { y^{1-s} \over 1-2s } \Gamma (s-1-m,y)
\eea
Since we are interested in solutions with exponential decay, we have set the homogeneous part of the solution to zero. Using the recursion relation for the incomplete $\Gamma$-function, 
\bea
\Gamma (a+1,y) = a \Gamma (a,y) + y^a \, e^{-y}
\eea
we may recast the result solely in terms of elementary functions and the exponential integral ${\rm Ei}_1 (y)= \Gamma (0,y)$. Simplifying the $y^s$ term, we find,
\bea
{y^s \over 2s-1} \Gamma (-s-m,y) = \gamma ^0 \, y^s \, {\rm Ei(y)} + \sum _{k=1}^{s+m} \gamma _k \, y^{s-k}  \, e^{-y} 
\eea
for rational coefficients $\gamma_0, \gamma_k$ whose value will not concern us here. Simplifying the $y^{1-s}$ term, the cases $s\geq m+2$ and $s \leq m+1$ must be distinguished, and we have,
\bea
s  \leq m+1 & & {y^{1-s} \over 1-2s} \, \Gamma (s-1-m,y)  = \gamma ^1 y^{1-s} \, {\rm Ei}(y) 
+ \sum _{k=1}^{1-s+m} \tilde \gamma _k \, y^{1-s+k} \, e^{-y}
\no \\
s  \geq m+1 & \hskip 0.5in &  {y^{1-s} \over 1-2s} \, \Gamma (s-1-m,y)  = 
 \sum _{k=0}^{s-m-2} \tilde \gamma _k \, y^{1-s+k} \, e^{-y}
\eea
Thus, the solution $\f_{s,m}(y)$ takes the form, 
\bea
\f_{s,m}(y) = \gamma ^0 \, y^s \, {\rm Ei}(y) +   \gamma ^1 y^{1-s} \, {\rm Ei}(y) 
+ \sum _{k =k_-} ^{s-1} \gamma _k \, y^k \, e^{-y} 
\eea
where $k_- = \min(-m,1-s)$ and the constant $\gamma ^1$ vanishes whenever $s \geq m+1$. Summing now over the range of $m$ in (5.8), and taking into account that $s \leq w-2$, we obtain, 
\bea
f_s(y) = f_+ \, y^s \, {\rm Ei}(y)+ f_- \, y^{1-s}\, {\rm Ei}(y) + 
\sum _{k =2-w} ^{w-3} f _k \, y^k \, e^{-y} 
\eea
When $s=w-2$ then we have $m \leq s$ and we must have $f_-=0$. In other cases, $f_\pm$ will generically be non-zero. This completes the proof of Theorem 1.3.

\subsection{The example of $C_{2,1,1}$}

We may work things out explicitly for the simplest case where the right side of the differential equation has a term which is non-linear in the Eisenstein series, and satisfies,
\bea
(\Delta -2) C_{2,1,1}= 9 E_4 - E_2^2
\eea
The reduced differential equation for the exponential part of the constant Fourier mode is,
\bea
\Big (\tau_2^2 \p_{\tau_2}^2  - 2 \Big ) \cC^{(e)} _{2,1,1} (\tau_2)  
=  
- { 4 \over 3} \sum_{n=1}^\infty  n^2 \sigma _{-3}(n)^2 
 e^{-4 \pi \tau_2 n} P_2(4 \pi \tau_2 n)^2 
\eea
where the polynomials $P_n$ were defined in (\ref{Pw}). Collecting the solution to the homogeneous part, and then summing over the particular solutions for each order in $N$ of the inhomogeneous part we find,
\bea
\cC^{(e)} _{2,1,1} (\tau_2) =
- { 4 \over 3} \sum_{n=1}^\infty  n ^2 \sigma _{-3}(n)^2 \f (4 \pi \tau_2 n)
\eea
where the function $\f$ satisfies the $n$-independent differential equation, 
\bea 
\Big (y^2 \p_y^2  - 2 \Big ) \f (y) =   e^{-y} P_2(y)^2
\eea
whose solution is as follows,
\bea
\f(y) =  { e^{-y} \over  y^2} 
\eea
This result is consistent with the result of Theorem \ref{t3} for $s=2$.

\section{The conjectured decomposition formula}
\setcounter{equation}{0}
\label{sec:7}

Theorem \ref{t2} states  that  the coefficient $c_{2-w}$ is given by a linear combination of products of pairs of odd $\zeta$-values whose weights add up to $2w-2$, but the theorem does not provide the coefficients $\gamma _\ell$ in (\ref{thm2}). In this section we shall obtain an explicit formula for  $\gamma _\ell$ in the decomposition of $c_{2-w}$. Using Maple, we have accumulated extensive numerical evidence for the validity of this formula, but we have no complete analytical proof. The decomposition formula reproduces all the special cases obtained in earlier work and  given in (\ref{spec}). 
To develop the required decomposition formula, we begin by proving the following Lemma.

\medskip
\noindent
\begin{lem}\label{l2}
The  linear combinations of  depth-two $\zeta$-values, 
\bea
\label{L2a}
\cS(M,N) & = & \zeta (2M-1, 2N+1) + \sum _{\ell=0}^{N-1} { \f_\ell (M,N) \, \Gamma (2M+2\ell) \over (2\ell+2)! \, \Gamma (2M-1)} 
\zeta (2M+2\ell, 2N-2\ell)
\no \\
\cT(M,N) & = & \zeta (2N+1,2M-1) + \sum _{\ell=0}^{N-1}  {  \f_\ell (M,N) \, \Gamma (2M+2\ell) \over (2\ell+2)! \, \Gamma (2M-1)}  \zeta (2N-2\ell,2M+2\ell)
\quad
\eea
obey the following relations.
\begin{description}
\itemsep=0in
\item (a)  For $N=0$, the function $\cS(M,0)$ is given by,
\bea
\label{L2b}
\cS(M,0)  =   \half(2M-1)\zeta (2M) - \half \sum_{j=1}^{2M-3} \zeta (j+1) \zeta (2M-1-j)
\eea 
\item (b) For $N \geq 1$, the functions  $\cS(M,N)$ and $\cT(M,N)$  are related as follows,
\bea
\label{L2e}
\cS(M,N)  & = & - \cT(M,N) + 
\zeta(2M-1) \zeta(2N+1)-\zeta(2M+2N)
\\ &&
+ \sum _{\ell=0}^{N-1}  { \f_\ell (M,N)  \, \Gamma (2M+2\ell) \over (2\ell+2)! \, \Gamma (2M-1)}  
\Big ( \zeta(2M+2\ell) \zeta(2N-2\ell) - \zeta (2M+2N) \Big)
\no
\eea
\item (c) The functions $\cS(M,N)$ and $\cT(M,N)$  reduce to a linear combination of products of two odd $\zeta$-values, plus  a  term valued in  $\pi^{2w-2} \QQ$, provided the coefficients $\f_\ell (M,N) $ are  given  in terms of the Euler polynomial $E_{2\ell+1}$  by,
\bea
\label{fi}
 \f_\ell (M,N)= - (2\ell+2) E_{2\ell+1}(0)
 \eea
 The coefficients $\f_\ell$ are integer-valued.
\item (d) For $N \geq 1$, the function $\cT(M,N)$ then evaluates as follows, 
\bea
\label{L2d}
\cT(M,N) & = & \sum _{\alpha =2} ^{2M+2N-2} \! \! (-)^{\alpha +1} \zeta (\alpha ) \zeta (2M+2N-\alpha)
\no \\ && \qquad \times
\sum _{n=0}^{2N-1} \half E_n(0) \binom{\alpha -1} { 2N-n} \binom{2M+n-2}{n}
\quad
\eea
\end{description}
\end{lem}
While there have been many investigations into the decomposition of linear combinations of double $\zeta$-values onto products of single $\zeta$-values, such as for example in \cite{GKZ, Mahide1, Mahide2}, the authors have not been able to find the precise statements of Lemma \ref{l2} in the literature.\footnote{We thank the referee for pointing out that the fact that the sum is totally reducible to a polynomial in the odd weight zeta values shows that the coefficient is a single-valued multiple zeta value.}

\subsection{Proof of Lemma \ref{l2}}

The proof of part {\sl (a)}  follows from a well-known formula of Euler, valid for integer $s \geq 2$,
\bea
\zeta (s,1) = {s \over 2} \zeta (s+1) - \half \sum _{j=1}^{s-2} \zeta (j+1) \zeta (s-j)
\eea

\sm

The proof of part {\sl (b)} is obtained by adding the expressions for $\cS(M,N)$ and $\cT(M,N)$ of (\ref{L2a}) and using Euler's reflection formula, valid for $s,t \in \CC$,
\bea
\label{4d1}
\zeta (s,t) + \zeta (t,s) & = & \zeta (s) \zeta (t) - \zeta (s+t)
\eea

\sm

The proof of part {\sl (c)} is more involved. Since the sum $\cS(M,N)+\cT(M,N)$ reduces to the product of two odd $\zeta$-values plus a term valued in  $\pi^{2w-2} \QQ$ by part {\sl (b)}, it will suffice to prove reducibility for $\cT(M,N)$.  To do so  we express the depth-two $\zeta$-functions in the definition of $\cT(M,N)$ in (\ref{L2a}) in terms of a double sum using (\ref{zeta2}),
\bea
\label{2f3}
\cT(M,N) = \half \sum _{m,n=1}^\infty F(m,n)
\eea
The function $F(m,n)$ is chosen to be symmetric in $m,n$, and  is given by,
\bea
\label{2f4}
F(m,n) & = & { 1 \over (m+n)^{2N+1} \, m^{2M-1} } +  { 1 \over (m+n)^{2N+1} \, n^{2M-1} } 
\\ &&
+ \sum _{\ell =0}^{N-1}  {  \f_\ell (M,N) \,  \Gamma (2M+2\ell) \over (2\ell+2)! \, \Gamma (2M-1)} 
\left ( { 1 \over (m+n)^{2N-2\ell} \, m^{2M+2\ell} } +  { 1 \over (m+n)^{2N-2\ell} \, n^{2M+2\ell} } \right )
\no \eea
Next, we obtain the conditions on the coefficients $\f_\ell (M,N)$ required for reducibility of $\cT(M,N)$, and show that  $\f_\ell(M,N)$ is  independent of $M,N$. To do so, we proceed by partial fraction decomposition of $F$ in the variable $m$ (or equivalently in $n$). Homogeneity of $F$ in the variables $m,n$ of degree $2M+2N$ restricts  the decomposition  to the following form,  
\bea
\label{2f6}
F(m,n) = \sum _{\alpha =0}^{2N} { \tilde f_\alpha (M,N) \over (m+n)^{2N+1-\alpha} \, n^{2M+\alpha-1}}
+ \sum _{\alpha=1}^{2M+2N-2} { f_\alpha (M,N) \over m^\alpha \, n^{2M+2N-\alpha}}
\eea
We determine the coefficients  $\f_\ell(M,N)$ by requiring $\tilde f_\alpha(M,N)=0$ for all $0 \leq \alpha \leq 2N$, namely the absence of double $\zeta$-values in the sum.  Although the system appears overdetermined with $2N+1$ conditions on $N$ variables $\f_\ell$ with $0 \leq \ell \leq N-1$, the symmetry in $m,n$ guarantees that it is uniquely solvable. To set $\tilde f_\alpha (M,N)=0$, we require that $(m+n)^{2N+1}F(m,n)$ and its $2N$ derivatives in $m$ vanish at $m=-n$. The vanishing of $\tilde f_{0}(M,N)$ is automatic. For $1 \leq \alpha \leq 2N$, the condition $\tilde f_\alpha=0$ is given by,
\bea
\sum _{\ell=0}^{N-1} { \f_\ell(M,N) \over 2 \ell +2} \left ( \delta _{2 \ell+1, \alpha } +
\left ( \bma \alpha \cr 2\ell+1 \cr \ema \right) \right ) =1
\eea
To solve these conditions, we proceed as follows. 
The sum over $\ell $ receives no contributions from $2\ell +1 > \alpha$. Therefore, the upper range of the sums may  be replaced by $[\alpha/2]$, and the conditions on $\f_\ell(M,N)$ are now independent of $M$ and $N$, and so are their solutions. Separating the equations for even and odd $\alpha$, we set respectively $\alpha = 2p+1$ and $\alpha =2p+2$. The  relations determining $\f_\ell= \f_\ell(M,N)$ are thus given as follows, 
\bea
\sum _{\ell=0}^p {  \f_\ell \over 2 \ell +2} \left ( \delta _{ \ell, p} +
\left ( \bma 2p +1 \cr 2\ell+1 \cr \ema \right) \right ) =
\sum _{\ell=0}^p { \f_\ell \over 2 \ell +2} 
\left ( \bma 2p+2 \cr 2\ell+1 \cr \ema \right) = 1 
\eea
These equations are equivalent to one another. We solve the second equation by multiplying by $x^{2p+1} /\Gamma (2p+2)$ and summing over all $p \geq 0$. The sum over $p$ may be carried out explicitly in terms of hyperbolic functions. Introducing the generating function,
\bea
\f(x) = \sum _{\ell=0}^ \infty { \f_\ell \, x^{2\ell} \over \Gamma (2 \ell+3)}
\eea
it is found that $\f(x)$ must satisfy $\f(0)=\half$ and the following linear differential equation, 
\bea
1 = x \f'(x) + \left ( 1 + { x \over \th (x)} \right ) \f(x)
\eea
Its integral is given by $\f(x) = ( e^x -1) x^{-1} (e^x +1)^{-1}$.
In terms of the Euler polynomials $E_n(t)$ which are defined by, 
\bea
{ 2 \, e^{xt} \over e^x +1} = \sum _{n=0}^\infty E_n(t) { x^n \over n!}
\eea
we find the explicit expression $  \f_\ell = - (2 \ell+2) E_{2\ell+1} (0)$. This completes the proof of part {\sl (c)}.
While $E_{2\ell+1}(0)$ are generally rational numbers, the coefficients $\f_\ell $ are integers. 

\sm

To prove part {\sl (d)}, we determine the coefficients $f_\alpha (M,N)$ in (\ref{2f6}). Using the expression we have obtained for $\f_\ell$ in (\ref{fi}), we readily find  $f_1(M,N)=0$. Furthermore, symmetry of $F$ in $m,n$ implies  the reflection property, $f_{2M+2N-\alpha} (M,N) = f_\alpha (M,N)$. Combining (\ref{2f3}) and (\ref{2f6}) with $\tilde f_\alpha (M,N)=0$, we obtain $\cT(M,N)$ by summing over $m,n$, 
\bea
\label{T1}
\cT(M,N) =  \sum_{\alpha=2}^{2M+2N-2} \half f_\alpha (M,N) \, \zeta (\alpha) \, \zeta (2M+2N-\alpha)
\eea
To compute the expressions for the remaining coefficients $f_\alpha(M,N)$, we match the poles in $n$ at $n=0$  between the expressions (\ref{2f4}) and (\ref{2f6}) with $\tilde f_\alpha (M,N)=0$, and we obtain, 
\bea
f_\alpha (M,N) = 
(-)^{\alpha+1}  \sum _{n=0}^{2N-1} E_{n}(0) \left ( \bma  \alpha -1 \cr 2N-n \cr \ema \right ) \left ( \bma 2M+n-2 \cr n \cr \ema \right ) 
\eea
Note that $E_0(0)=1$, and $E_{2\ell}(0)=0$ for $\ell \in \NN$, which simplifies the above sum to ranging only over $n=0$ and  the odd positive integers. Substituting the expression for $f_\alpha(M,N)$ into (\ref{T1}) gives  (\ref{L2d})  and thereby proves part (d).

\subsection{Derivation of the decomposition formula}

The starting point for the decomposition formula is the expression for $c_{2-w}$ in terms of double $\zeta$-values, given in (\ref{C2W}), (\ref{C2W0}) and (\ref{ZZ}) of Theorem \ref{t1}. We change summation variable
in the definition of $Z(a_1, a_2, a_3)$ in (\ref{ZZ}) from $\ell $ to $m=k+\ell-1$, so that we have,
\bea
Z(a_1,a_2,a_3)
& = & 
\sum_{m=1}^{2a_1-1} \sum _{k=k_-}^{k_+} 
\binom{a_1+a_2-k-1}{a_2-1} \binom{a_1+a_2-m+k-2}{a_2-1} 
\no \\ && \hskip 0.4in \times 
\binom{m-1}{k-1} \binom{2w-m-3}{w-k-1} \, \zeta (2w-m-2, m)
\eea
where $k_\pm$ are defined by $k_+ = \min(a_1,m)$ and $k_- = \max(1,m-a_1+1)$. Next, we eliminate 
the contributions  of the odd-odd double $\zeta$-values, which arise when $m$ is odd,  in favor of the function $\cS$ plus even-even $\zeta$-values. The most interesting formula is obtained upon symmetrizing in $a_2$ and $a_3$, and we find, 
\bea
\label{L2c}
Z(a_1, a_2, a_3) + Z(a_1, a_3, a_2) & = & \sum _{n=0}^{a_1-1} \Big ( Z_n (a_1, a_2, a_3) 
+ Z_n (a_1, a_3, a_2) \Big ) \cS(w-1-n,n)
\no \\ &&
+ \sum _{n=1}^{a_1-1} X_n (a_1, a_2, a_3) \zeta(2w-2n-2,2n)
\eea
where $Z_\alpha (a_1, a_2, a_3)$ is given in (\ref{sumL4}), and the function $X_n (a_1, a_2, a_3)$ is defined by, 
\bea
X_n(a_1, a_2, a_3) & = & 
\sum _{\ell =n}^{a_1-1} E_{2\ell-2n+1}(0) Z_\ell(a_1, a_2, a_3) \binom{2w-2n-3}{2w-2\ell-4}
\no \\ &&
+\sum _{k=k_-'}^{k_+'} 
\binom{a_1+a_2-k-1}{a_2-1} \binom{a_1+a_2-2n+k-2}{a_2-1} 
\no \\ && \hskip 0.4in \times 
\binom{2n-1}{k-1} \binom{2w-2n-3}{w-k-1} + (a_2 \leftrightarrow a_3)
\eea
where $k_\pm'$ are defined by $k_+' = \min(a_1,2n)$ and $k_-' = \max(1,2n-a_1+1)$. The first sum on the right arises from the elimination in (\ref{ZZ}) of the odd-odd double $\zeta$-values in favor of $\cS$, while the second term arises from the contribution of the even-even $\zeta$-values to (\ref{ZZ}). The symmetrization in $a_2, a_3$ applies to the entire expression. The formula for $X_n$ may be rendered more explicit by expressing  $Z_\ell$ as the sum given in (\ref{sumL4}), changing variables $n \to a_1 -n$ and $\ell \to a_1 -\ell$, and recognizing that the second sum may be combined with the first corresponding to the Euler polynomial with index zero, $E_0(0)=1$, and we find, 
\bea
X_{a_1-n} (a_1, a_2, a_3) 
& = & 
\sum _{k=0}^{a_1-1} \sum _{\ell=0}^{a_1-1} \theta (2n-k-\ell-1) E_{2n-k-\ell-1} (0)
\binom{a_2-1+k}{k} \binom{a_2-1+\ell}{\ell} 
\no \\ && \times 
\binom{2a_1-k-\ell-2}{a_1-k-1} 
\binom{2a_+2a_3+k+\ell-2}{a_2+a_3+k-1} 
\binom{2a_2+2a_3+2n-3}{2n-k-\ell-1}
\no \\ &&  
+ (a_2 \leftrightarrow a_3)
\eea
Note that the index of the Euler polynomial is now allowed to run over even and odd integers though, of course, we have $E_{2k}(0)=1$ for all $k\geq 1$.
 
\medskip 
\noindent
\begin{conj} \label{conj2}
For $a_1, a_2, a_3  \in \NN$ and for $1 \leq n \leq a_1-1$ the following identities hold,  
\bea
X_n(a_1, a_2, a_3)=0 
\eea
\end{conj}

\sm

Using MAPLE,  the conjecture has been verified by showing that $X_n(a_1, a_2, a_3)=0$ as functions of $a_2$ and $a_3$  for all $n$ in the range $1 \leq n \leq a_1-1$ and for   $1 \leq a_1 \leq 65$.
An analytical proof of the conjecture is, however, outstanding.

\sm

Assuming the validity of Conjecture \ref{conj2}, we prove the {\sl decomposition formula} of Conjecture \ref{conj1} by expressing $c_{2-w}$ in terms of $\cS$. To this end, we combine the formulas (\ref{C2W}), (\ref{C2W0}), (\ref{L2c}) and the result of Conjecture \ref{conj2}, to obtain, 
\bea
c_{2-w} = c_{2-w}^0 \zeta (2w-2) + 2 \sum _{\sigma \in \mS_3} \sum _{n=0}^{a_1-1} Z_n(a_1, a_2, a_3) \cS(w-1-n,n)
\eea
Evaluating $\cS$ with the help of formulas (\ref{L2b}), (\ref{L2e}) and (\ref{L2d}) of Lemma \ref{l2} allows us to express $\cS$ in terms of a sum over pairs of odd $\zeta$-values whose weights add up to $2w-2$, as well as terms taking values in $\pi^{2w-2} \QQ$ (arising from the contributions of even $\zeta$-values in (\ref{L2b}), (\ref{L2e}) and (\ref{L2d}), the latter for even $\alpha$). By Theorem \ref{t2}, the contributions valued in $\pi^{2w-2} \QQ$ cancel. The remaining contributions give $\gamma _k$ of (\ref{sumL3}) as follows. The first term on the right in (\ref{sumL3}) arises from $\cS(w-1,0)$ given by (\ref{L2b}), the second term from the second term on the right in the sum of $\cS$ and $\cT$ in (\ref{L2e}), and the remaining terms from the terms for odd values of $\alpha$ in the decomposition of $\cT$ in (\ref{L2d}). This completes the proof of the decomposition formula of Conjecture \ref{conj1}) assuming the validity of Conjecture \ref{conj2}.

\medskip

\subsection*{Acknowledgments}

It is a pleasure to thank Michael Green and Pierre Vanhove for many discussions and earlier collaborations  out of which the present work has grown. We are also happy to thank Justin Kaidi for some assistance in articulating Lemma 1 and for helpful comments on the manuscript.  ED is grateful to the Kavli Institute for Theoretical Physics in Santa Barbara for their hospitality while part of this work was being completed. This research  is supported in part by the National Science Foundation under research grants PHY-16-19926 and PHY-1125915 for ED and DMS 1701638 for WD. ED is pleased to acknowledge support by a Fellowship from the Simons Foundation.  

\medskip


\end{document}